\documentclass[letterpaper]{article} 
\usepackage{aaai2026}  
\usepackage{times}  
\usepackage{helvet}  
\usepackage{courier}  
\usepackage[hyphens]{url}  
\usepackage{graphicx} 
\usepackage{natbib}  
\usepackage{caption} 
\usepackage{algorithm}
\usepackage{algorithmic}
\usepackage{amsthm}
\usepackage{amsmath}
\usepackage{amssymb}
\usepackage{enumitem}
\usepackage{atnguyen_new}
\usepackage{soul}
\usepackage{mathtools}
\usepackage{newfloat}
\usepackage{listings}
\DeclareCaptionStyle{ruled}{labelfont=normalfont,labelsep=colon,strut=off} 
\floatstyle{ruled}
\newfloat{listing}{tb}{lst}{}
\floatname{listing}{Listing}
\newcommand{\LineComment}[1]{\STATE \textit{\textcolor{gray}{/* #1 */}}}

\title{Provable Data-driven Projection Method for Quadratic Programming}
\author {
    Anh Tuan Nguyen\textsuperscript{\rm 1}\thanks{work done while visiting TTIC},
    Viet Anh Nguyen\textsuperscript{\rm 2},
}
\affiliations {
    \textsuperscript{\rm 1}Carnegie Mellon University, 
    \textsuperscript{\rm 2}Chinese University of Hong Kong\\
    atnguyen@cs.cmu.edu, nguyen@se.cuhk.edu.hk
}

\usepackage{bibentry}

\begin{document}

\maketitle

\begin{abstract}
    Projection methods aim to reduce the dimensionality of the optimization instance, thereby improving the scalability of high-dimensional problems. Recently, \citet{sakaue2024generalization} proposed a data-driven approach for linear programs (LPs), where the projection matrix is learned from observed problem instances drawn from an application-specific distribution of problems. We analyze the generalization guarantee for the data-driven projection matrix learning for convex quadratic programs (QPs). Unlike in LPs, the optimal solutions of convex QPs are not confined to the vertices of the feasible polyhedron, and this complicates the analysis of the optimal value function. To overcome this challenge, we demonstrate that the solutions of convex QPs can be localized within a feasible region corresponding to a special active set, utilizing Carath\'{e}odory's theorem. Building on such observation, we propose the \textit{unrolled active set method}, which models the computation of the optimal value as a Goldberg-Jerrum algorithm with bounded complexities, thereby establishing learning guarantees. We then extend our analysis further to other settings, including learning to match the optimal solution and an input-aware setting, where we learn to map QP problem instances to projection matrices.
\end{abstract}

\begin{links}
    \link{Extended version}{https://arxiv.org/pdf/2509.04524}
\end{links}

\section{Introduction}
Linear programs (LPs) and the more general quadratic programs (QPs) are simple forms of convex optimization problems, yet they play crucial roles in many industrial~\cite{gass2003linear, dostal2009optimal} and scientific domains~\cite{amos2202tutorial}. Practical LP and QP instances are often computationally intensive to solve due to their enormous problem size, which can encompass millions of variables and constraints. As a result, accelerating solving approaches for large-scale LPs and QPs are important directions in the operations research literature. Two of the most prominent approaches include accelerated solvers and dimensionality reduction methods. Accelerated solvers focus on improving the speed of widely used solvers on large-scale problems via parallelization, randomization, or wisely leveraging cheap first-order (i.e., gradient) information, among other approaches. Some of the recent advances include parallelized simplex methods \cite{huangfu2018parallelizing}, randomized interior-point methods \cite{chowdhury2022faster}, and primal-dual hybrid gradient methods \cite{applegate2021practical}. 

Another complementary, solver-agnostic approach for large-scale LPs and QPs is the dimensionality reduction technique, which generalizes the idea of reducing the size of problem instances while preserving the properties of the objective values and variables. A promising approach is through random projections \cite{d2020random,vu2019random,vu2018random}, where a random projection matrix is used to map the variables and feasible regions of the original problem instances into a low-dimensional space, forming projected problem instances that can be solved much faster. The solutions of projected problem instances can then be mapped back to the original space, in the hope that their quality is sufficiently comparable to that of the optimal solution of the original problem instances. Importantly, this solver-agnostic approach can be combined with accelerated solvers to further improve the solving of large-scale LPs and QPs.

However, random projection matrices overlook the geometric properties of the problem instances, and this oversight potentially leads to inferior solution quality of the projected problem instances compared to that of the original problem instances. Recently, \citet{sakaue2024generalization} proposed a data-driven approach for learning the projection matrix, specifically targeting LPs. Assume that there are not one, but multiple LPs $\boldsymbol{\pi}_\textup{LP} = (\boldsymbol{c}, \boldsymbol{A}, \boldsymbol{b}) \in \Pi_\textup{LP} \subset \bbR^{n} \times \bbR^{n \times m} \times \bbR^m$ that have to be solved in the form 
\[
    \textup{OPT}(\boldsymbol{\pi}_\textup{LP}) = \min_{\boldsymbol{x} \in \bbR^n} \{\boldsymbol{c}^\top \boldsymbol{x} \mid \boldsymbol{Ax} \leq \boldsymbol{b} \}.
\]
The parameters $\boldsymbol{\pi}_\textup{LP}$ are drawn from some application-specific and potentially unknown problem distribution $\cD_\textup{LP}$ over $\Pi_\textup{LP}$. \citet{sakaue2024generalization} proposed to learn the projection matrix $\boldsymbol{P} \in \cP \subset \bbR^{n \times k}$, where $k \ll n$ is the dimensionality of the projection space, by minimizing the expected optimal objective of the projected LPs $\bbE_{\boldsymbol{\pi}_\textup{LP} \sim \cD_\textup{LP}} [\ell_\textup{LP}(\boldsymbol{P}, \boldsymbol{\pi}_\textup{LP})]$, where
\[
    \ell_\textup{LP}(\boldsymbol{P}, \boldsymbol{\pi}_\textup{LP}) = \min_{\boldsymbol{y} \in \bbR^k} \{\boldsymbol{c}^\top \boldsymbol{P}\boldsymbol{y}\mid \boldsymbol{A}\boldsymbol{Py}\leq \boldsymbol{b} \}
\]
is the optimal objective value of the projected LP. Because $\cD$ is unknown, minimizing $\bbE_{\boldsymbol{\pi}_\textup{LP} \sim \cD_\textup{LP}} [\ell_\textup{LP}(\boldsymbol{P}, \boldsymbol{\pi}_\textup{LP})]$ is intractable, and instead we learn $\boldsymbol{P}$ through \textit{empirical risk minimization (ERM)} using LP problem instances drawn from $\cD_\textup{LP}$. It is easy to see that $\ell_\textup{LP}(\boldsymbol{P}, \boldsymbol{\pi}_\textup{LP})$ upper-bounds $\textup{OPT}(\boldsymbol{\pi}_\textup{LP})$, and therefore the smaller $\bbE_{\boldsymbol{\pi}_\textup{LP} \sim \cD_\textup{LP}} [\ell_\textup{LP}(\boldsymbol{P}, \boldsymbol{\pi}_\textup{LP})]$ is, the closer the quality of solutions of projected problem instances to that of original problem instances. Along with promising empirical results, \citet{sakaue2024generalization} provided generalization guarantees for learning $\boldsymbol{P}$ via ERM by analyzing the learning-theoretic complexity (i.e., pseudo-dimension \cite{pollard1984convergence}) of the corresponding loss function class $\cL_\textup{LP} = \{\ell_{\boldsymbol{P}}:\Pi_\textup{LP} \rightarrow [-H, 0] \mid \boldsymbol{P} \in \cP\}$, where $\ell_{\boldsymbol{P}}(\boldsymbol{\pi}_\textup{LP}) \coloneq \ell_\textup{LP}(\boldsymbol{P}, \boldsymbol{\pi}_{\textup{LP}})$, and $H$ is some real-valued upper-bound for the function class.

Inspired by this success, a natural direction is to extend this framework to convex QPs. Similarly, given QP problem instances $\boldsymbol{\pi} = (\boldsymbol{Q}, \boldsymbol{c}, \boldsymbol{A}, \boldsymbol{b}) \in \Pi \subset \bbR^{n \times n} \times \bbR^n \times \bbR^{m \times n} \times \bbR^n$ coming from an application-specific, unknown problem distribution $\cD$ over $\Pi$, the idea is to learn a projection matrix $\boldsymbol{P} \in \cP \subset \bbR^{n \times k}$ with $k \ll n$ that achieves small population loss $\bbE_{\boldsymbol{\pi} \sim \cD} [\ell(\boldsymbol{P}, \boldsymbol{\pi})]$ via ERM, where
\[
    \ell(\boldsymbol{P}, \boldsymbol{\pi}) = \min_{\boldsymbol{y} \in \bbR^k} \biggr\{ \frac{1}{2}\boldsymbol{y}^\top \boldsymbol{P}^\top \boldsymbol{Q}\boldsymbol{P\boldsymbol{y}} + \boldsymbol{c}^\top \boldsymbol{Py} \mid \boldsymbol{A}\boldsymbol{Py} \leq \boldsymbol{b} \biggr\}.
\]
Again, to ensure the generalization guarantee for $\boldsymbol{P}$ learned via ERM, we need to analyze the function class $\cL = \{\ell_{\boldsymbol{P}}: \Pi \rightarrow [-H, 0] \mid \boldsymbol{P} \in \cP\}$, where $\ell_{\boldsymbol{P}}(\boldsymbol{\pi}) \coloneqq \ell(\boldsymbol{P}, \boldsymbol{\pi})$.

At first glance, the extension to QPs may seem straightforward because the previous ideas seem readily applicable to the form of QPs, and the gradient update can also be derived using the envelope theorem. However, the optimal solutions of QPs exhibit fundamentally different geometrical structures, and it turns out that extending the existing theoretical framework to QPs requires developing new tools tailored to these specific problems. 

\paragraph{Contributions.} We formalize the data-driven projection method for convex QPs and analyze the generalization guarantees of learning the projection matrix. Our contributions can be summarized as follows:
\begin{enumerate}
    \item We establish generalization guarantees for the data-driven learning projection matrix for QPs in~Theorem~\ref{thm:pdim-original-bound}. Our new result is more general and strictly tighter than the previous bound proposed by~\citet{sakaue2024generalization}, which applies only to LPs. For completeness, we also instantiate a lower bound for the convex QPs case in Proposition~\ref{prop:lower-bound}.
    \item We propose and analyze a novel learning scenario, where the goal is to match the optimal solution in Section~\ref{sec:match-solution}. This setting is particularly useful in practical applications where the focus is on the solution to be implemented. The guarantee result is presented in Theorem~\ref{thm:learning-to-match-solution}.
    \item We consider input-aware settings, where we learn a neural network that maps a convex QP to a tailored projection matrix, as described in Section~\ref{sec:input-aware}. The guarantee result is presented in Theorem \ref{thm:input-aware-setting}.
\end{enumerate}

\textbf{Technical challenges and overviews.} For any LP with parameters $\boldsymbol{\pi}_\textup{LP} = (\boldsymbol{c}, \boldsymbol{A}, \boldsymbol{b})$ and any projection matrix $\boldsymbol{P}$, the solution of the projected LP is one of the vertices of the feasible polyhedron. Leveraging such observation, \citet{sakaue2024generalization} describes the computation of the projected LP's optimal value $\ell_\textup{LP}(\boldsymbol{P}, \boldsymbol{\pi}_\textup{LP})$ by enumerating all potential vertices, and identifies the vertex $\boldsymbol{y}^*$ that yields the lowest objective $\boldsymbol{c}^\top \boldsymbol{P}\boldsymbol{y}^*$. The computation of $\ell_\textup{LP}(\boldsymbol{P}, \boldsymbol{\pi}_\textup{LP})$ can then be described by a bounded number of distinct conditional statements involving polynomials in the entries of $\boldsymbol{P}$; see Section~\ref{sec:GJ} for details.

This favorable property, however, does not extend to QPs, as the solution of QPs can be anywhere within the feasible polyhedron, not just at its vertices. This makes directly locating the solution and calculating the optimal objective $\ell(\boldsymbol{P}, \boldsymbol{\pi})$ very challenging. To overcome this issue, we propose a four-step analytical approach. First, we will construct a perturbed objective $\ell_{\boldsymbol{P}, \gamma}(\boldsymbol{\pi})$ that is well-behaved and can approximate $\ell(\boldsymbol{P}, \boldsymbol{\pi})$ with arbitrarily precision (Lemma \ref{lm:perturbed-QP-approximating-QP}, Proposition \ref{prop:perturbed-OQP-favourable-structure}). Second, we leverage the structure of this perturbed problem to develop the \textit{unrolled active set method}, an algorithm that exactly computes its optimal value (Lemma \ref{lm:unrolled-active-method-correctness}). Third, we demonstrate that our method can be framed as a GJ algorithm with bounded complexities (Lemma \ref{lm:gj-complexity-unrolled-active-set}), which enables us to bound the pseudo-dimension of the perturbed function class. Finally, by relating the perturbed objective to the original, we extend this bound to the original QP loss function, thereby proving our main generalization guarantee (Theorem \ref{thm:pdim-original-bound}).

\section{Related Works}
\paragraph{Projection methods for LPs and QPs.} Projection methods aim to accelerate the solution of LPs and QPs by reducing the size of the problem instances. Previous works have investigated random projection to reduce the number of constraints~\cite{vu2019random, poirion2023random} and variables~\cite{akchen2025column}. Recently, \citet{sakaue2024generalization} and \citet{iwata2025learning} considered a data-driven approach, learning the projection matrix for a specific distribution of problems instead of a random projection, specifically targeting LPs. Our paper extends this framework to convex QPs.

\paragraph{Learning to optimize.} Learning to optimize leverages machine learning to develop optimization methods, i.e., by predicting an initial solution for the exact algorithm, approximating the exact solution directly, or adapting specific components of optimization algorithms \cite{chen2022learning, amos2023tutorial, bengio2021machine}. Learning to project for LPs \cite{sakaue2024generalization, iwata2025learning} and convex QPs (this work) belongs to this broad category, where the learned projection matrices are used to accelerate off-the-shelf solvers and produce approximate solutions that are guaranteed to be feasible, unlike prior methods that approximate optimal solutions directly using neural networks.

\paragraph{Data-driven algorithm design.} Data-driven algorithm design \citep{balcan2020data, gupta2020data} is an emerging algorithm design paradigm that proposes adapting algorithms by configuring their hyperparameters or internal components to the specific set of problem instances they must solve, rather than considering the worst-case problem instances. Assuming that there is an application-specific, potentially unknown problem distribution from which the problem instances are drawn, data-driven algorithm design aims to maximize its empirical performance using the observed problem instances, with the hope that the adapted algorithm will perform well on future problem instances drawn from the same problem distribution. Data-driven algorithm design is an active research direction in both empirical validation and theoretical analysis across various domains, including sketching and low-rank approximation \citep{indyk2019learning, bartlett2022generalization, li2023learning}, (mixed) integer linear programming \citep{balcan2018learning, li2023learning}, tuning regularization hyperparameters \citep{balcan2022provably, balcan2023new}, and other general frameworks for theoretical analysis in data-driven settings \citet{bartlett2022generalization, balcan2025algorithm, balcan2025sample}. Data-driven projection methods for LPs \citep{sakaue2024generalization} and QPs are specific instances of data-driven algorithm design.

\section{Backgrounds on Learning Theory}
\subsection{Pseudo-dimension}
    We recall the notion of \textit{pseudo-dimension}, the primary learning-theoretic complexity measure of this work.

    \begin{definition}[Pseudo-dimension, \citealp{pollard1984convergence}] \label{def:pseudo-dimension}
        Consider a real-valued function class $\cL$, of which each function $\ell$ takes input $\boldsymbol{\pi}$ in $\Pi$ and output $\ell(\boldsymbol{\pi}) \in [-H, 0]$. Given a set of inputs $S = (\boldsymbol{\pi}_1, \dots, \boldsymbol{\pi}_N) \subset \Pi$, we say that $S$ is shattered by $\cL$ if there exists a set of real-valued threshold $r_1, \dots, r_N \in \bbR$ such that $\abs{\{(\sign(\ell(\boldsymbol{\pi}_1) - r_1), \dots, \sign(\ell(\boldsymbol{\pi}_N) - r_N)) \mid \ell \in \cL\}} = 2^N$. The pseudo-dimension of $\cL$, denoted as $\Pdim(\cL)$, is the maximum size $N$ of a input set that $\cL$ can shatter.
    \end{definition}

    It is widely known from the learning theory literature that if a real-valued function class has bounded pseudo-dimension, then it is PAC-learnable with ERM.
    
    \begin{theorem}[\citealp{pollard1984convergence}] \label{thm:pdim-to-generalization}
        Consider a real-valued function class $\cF$, of which each function $\cL$ takes input $\boldsymbol{\pi}$ in $\Pi$ and output $\ell(\boldsymbol{\pi}) \in [-H, 0]$. Assume that $\Pdim(\cL)$ is finite. Then given $\epsilon > 0$ and $\delta \in (0, 1)$, for any $N \geq M(\delta, \epsilon)$, where $M(\delta, \epsilon) = \cO\left(\frac{H^2}{\epsilon^2}(\Pdim(\cL) + \log(1/\delta))\right)$, with probability at least $1 - \delta$ over the draw of $S = (\boldsymbol{\pi}_1, \dots, \boldsymbol{\pi}_N) \sim \cD^N$, where $\cD$ is a distribution over $\Pi$, we have 
        \[
            \bbE_{\boldsymbol{\pi} \sim \cD}[\hat{\ell}_S(\boldsymbol{\pi})] \leq \inf_{\ell \in \cL} \bbE_{\boldsymbol{\pi} \sim \cD}[\ell(\boldsymbol{\pi})] + \epsilon.
        \]
        Here, $\hat{\ell}_S \in \arg \min_{\ell \in \cL}\frac{1}{N}\sum_{i = 1}^N \ell(\boldsymbol{\pi}_i)$ is the ERM minimizer.   
    \end{theorem}

\subsection{Goldberg-Jerrum Framework}  \label{sec:GJ}

    Goldberg-Jerrum (GJ) framework is originally proposed by \citet{goldberg1993bounding}, and is later refined by \citet{bartlett2022generalization}. It is a convenient framework to establish a pseudo-dimension upper-bound for parameterized function classes, of which the computation can be described by a \textit{GJ algorithm} using conditional statements, intermediate values, and outputs involving rational functions of their parameters. The formal definition of the GJ algorithm can be described as follows.
    
    \begin{definition}[GJ algorithm, \cite{bartlett2022generalization}]
        \label{def:GJ-algorithm}
         A \textbf{GJ algorithm} $\Gamma$ operates on real-valued inputs, and can perform two types of operations:
        \begin{itemize}
            \item Arithmetic operators of the form $v'' = v \odot v'$, where $\odot \in \{ +, -, \times, \div\}$, and
            \item Conditional statements of the form ``if  $v \geq 0 \ldots$ else $\ldots$".
        \end{itemize}
        In both cases, $v$ and $v'$ are either inputs or values previously computed by the algorithm.
    \end{definition}

    The immediate values $v, v', v''$ computed by the GJ algorithm are rational functions (fractions of two polynomials) of its parameters. The complexities of the GJ algorithm are measured by the highest degree of rational functions it computes and the number of distinct rational functions that appear in the conditional statements. The formal definition of its complexities is as follows.
        
    \begin{definition}[Complexities of GJ algorithm, \cite{bartlett2022generalization}]
        \label{definition:gj-algorithm-complexities}
         The \textbf{degree} of a GJ algorithm is the maximum degree of any rational function that it computes of the inputs. The \textbf{predicate complexity} of a GJ algorithm is the number of distinct rational functions that appear in its conditional statements. Here, the degree of rational function $f(\boldsymbol{x}) = \frac{g(\boldsymbol{x})}{h(\boldsymbol{x})}$, where $g$ and $h$ are two polynomials in $\boldsymbol{x}$, is $\deg(f) = \max\{\deg(g), \deg(h)\}$.
    \end{definition}

    The following theorem asserts that if any function class for which the function's computation can be described by a GJ algorithm with bounded degree and predicate complexities, then the pseudo-dimension of that function class is also bounded.
    
    \begin{theorem}[{\citet[theorem~3.3]{bartlett2022generalization}}] \label{thm:gj}
        Suppose that each function $\ell_{\boldsymbol{P}} \in \cL$ is specified by $n$ real parameters $\boldsymbol{P} \in \bbR^n$. Suppose that for every $\boldsymbol{\pi} \in \Pi$ and $r \in \bbR$, there is a GJ algorithm $\Gamma_{\boldsymbol{\pi}, r}$ that, given $\ell_{\boldsymbol{P}} \in \cL$, returns ``true" if $\ell_{\boldsymbol{P}}(\boldsymbol{\pi}) \geq r$ and ``false" otherwise. Assume that $\Gamma_{\boldsymbol{\pi}, r}$ has degree $\Delta$ and predicate complexity $\Lambda$. Then, $\Pdim(\cL) = \cO(n\log(\Delta\Lambda))$.
    \end{theorem}

        Note that the GJ algorithm $\Gamma_{\boldsymbol{\pi}, r}$ described above corresponds to each fixed input $\boldsymbol{\pi}$ and threshold value $r$. The input of the GJ algorithm $\Gamma_{\boldsymbol{\pi}, r}$ is the hyperparameters $\boldsymbol{P}$ (the projection matrix in our case) parameterizing $\ell_{\boldsymbol{P}}$, and the intermediate values and conditional statements involve in rational functions of $\boldsymbol{P}$. Moreover, the GJ framework only serves as a tool for analyzing \textit{learning-theoretic complexity} (e.g., pseudo-dimension) of the parameterized function class $\cL$, and does not describe how the function $\ell_{\boldsymbol{P}}(\boldsymbol{\pi})$ is computed in practice. In our framework, the computation of $\ell_{\boldsymbol{P}}(\pi)$ utilizes our proposed \textit{unrolled active set method} in Algorithm \ref{alg:unrolled-active-set-method}, which we show to be a GJ algorithm with bounded complexities in Lemma~\ref{lm:gj-complexity-unrolled-active-set}. In practice, it might be computed using the \textit{active set method} \cite{nocedal2006numerical} or interior-point method \cite{dikin1967iterative} for computational efficiency; however, these methods cannot be cast as GJ algorithms.
    
\section{Problem Settings} \label{sec:problem-setting}
This section formalizes the problem of learning the projection matrix for QPs in the data-driven setting. 
\subsection{Original QPs and Projected QPs}
Consider the \textit{original QPs} (OQPs) $\boldsymbol{\pi} = (\boldsymbol{Q}, \boldsymbol{c}, \boldsymbol{A}, \boldsymbol{b}) \in \Pi \subset \bbR^{n \times n} \times \bbR^n \times \bbR^{m \times n} \times \bbR^m$ with inequality constraints:
\[  \label{eq:OQP}
    \textup{OPT}(\boldsymbol{\pi}) = \min_{\boldsymbol{x \in \bbR^n}}\left\{\frac{1}{2}\boldsymbol{x}^\top \boldsymbol{Q} \boldsymbol{x} + \boldsymbol{c}^\top\boldsymbol{x} \mid \boldsymbol{Ax \leq b}\right\},
    \tag{OQP}
\]
where $\boldsymbol{Q}$ is a positive semi-definite (PSD) matrix, while $n$ and $m$ are the number of variables and constraints, respectively. Here, we assume that the variable size $n$ and the number of constraints $m$ are large, and solving the OQP is a computationally expensive task. The core idea of the projection method evolves around a \textit{full column-rank} projection matrix $\boldsymbol{P} \in \cP \subset \bbR^{n \times k}$, where $k \ll n$ is the projection dimension. Setting $\boldsymbol{x} = \boldsymbol{Py}$, we obtain the \textit{projected QPs} (PQPs) corresponding to the OQPs $\boldsymbol{\pi}$ and the projection matrix $\boldsymbol{P}$
\[ \label{eq:PQP}
    \ell(\boldsymbol{P}, \boldsymbol{\pi}) = \min_{\boldsymbol{y} \in \bbR^k} \left\{\frac{1}{2}\boldsymbol{y}^\top \boldsymbol{P}^\top \boldsymbol{Q}\boldsymbol{P}\boldsymbol{y} + \boldsymbol{c}^\top \boldsymbol{Py}\mid \boldsymbol{APy} \leq \boldsymbol{b}\right\}.
    \tag{PQP} 
\]
Similar to prior works \cite{sakaue2023improved, vu2019random}, we make the following assumptions for the OQPs.
\begin{assumption}[Regularity conditions] \label{asmp:regularity}
    The OQPs: 
    \begin{enumerate}[label = (\arabic*)]
    \item take inequality-constrained form as \eqref{eq:OQP}, 
    \item have $\boldsymbol{0}_n \in \bbR^n$ as a feasible point, 
    \item have the feasible region is bounded by $R$ in the sense that $\| \boldsymbol{x} \|_2 \le R$ for any feasible $\boldsymbol{x}$, and 
    \item have bounded optimal objective value from $[-H, 0]$, for some constant positive $H$.
    \end{enumerate}
\end{assumption}
\begin{remark}
    As discussed in previous work \cite{sakaue2024generalization, vu2019random}, Assumption \ref{asmp:regularity} is not restrictive. First, any QPs that also have equality assumptions can also be converted into the inequality form (Assumption \ref{asmp:regularity}.1 by considering the null space of the equality constraints (see Appendix C, \cite{sakaue2024generalization} for details). For Assumption \ref{asmp:regularity}.2, one can instead assume that there exists a feasible point $\boldsymbol{x}_0$, and linearly translate the feasible region so that $\boldsymbol{x}_0$ coincides with $\boldsymbol{0}_n$, without changing the form of QPs. Assumption \ref{asmp:regularity}.3 is standard from the optimization literature \cite{vu2019random}. Lastly, Assumption \ref{asmp:regularity}.4 is simply a consequence of Assumption \ref{asmp:regularity}.1 and \ref{asmp:regularity}.2.
\end{remark}
Under Assumption \ref{asmp:regularity}, the PQPs also have a favorable structure, which can be formalized as follows.
\begin{proposition} \label{prop:PQP-properties}
    Under Assumption \ref{asmp:regularity}, then for any OQP $\boldsymbol{\pi}$ and projection matrix $\boldsymbol{P}$, (i) the corresponding PQP has $\boldsymbol{0}_k$ as a feasible point, and (ii) $\ell(\boldsymbol{P}, \boldsymbol{\pi})$ is lower bounded by $\textup{OPT}(\boldsymbol{\pi})$, and therefore takes a value between $[-H, 0]$.
\end{proposition}
\proof
    Since $\boldsymbol{0}_n$ is a feasible point of OQP $\boldsymbol{\pi}$, $\boldsymbol{y} = \boldsymbol{0}_k$ satisfies $\boldsymbol{AP}\boldsymbol{y} \leq \boldsymbol{b}$, meaning that $\boldsymbol{0}_k$ is a feasible point of PQP. Moreover, let $\boldsymbol{y}^*$ be an optimal solution of PQP, then $\boldsymbol{x}' = \boldsymbol{P}\boldsymbol{y}^*$ is a feasible point of OQP, thus $\textup{OPT}(\boldsymbol{\pi}) \leq \ell(\boldsymbol{P}, \boldsymbol{\pi})$.
\qed

\subsection{Data-driven Learning of the Projection Matrix} \label{sec:data-driven-learning-projection-matrix}
In the data-driven setting, we assume that there is an application-specific and potentially unknown problem distribution $\cD$ over the set of QPs $\Pi$. The optimal projection matrix $\boldsymbol{P}$ minimizes the population PQPs' optimal objective value
\[
    \boldsymbol{P}^*_\cD \in \argmin_{\boldsymbol{P} \in \cP} \bbE_{\boldsymbol{\pi} \sim \cD} [\ell(\boldsymbol{P}, \boldsymbol{\pi})].
\]
From Proposition \ref{prop:PQP-properties}, we know that the smaller $\bbE_{\boldsymbol{\pi} \sim \cD} [\ell(\boldsymbol{P}, \boldsymbol{\pi})]$, the closer the optimal objective value of PQP $\ell(\boldsymbol{P}, \boldsymbol{\pi})$ to $\textup{OPT}(\boldsymbol{\pi})$, and the better $\boldsymbol{P}$. However, since $\cD$ is unknown, we instead learn $\boldsymbol{P}$ using the observed PQPs $S = \{\boldsymbol{\pi}_1, \dots, \boldsymbol{\pi}_N\}$ drawn i.i.d.~from $\cD$ via ERM
\[
    \hat{\boldsymbol{P}}_S \in \argmin_{\boldsymbol{P} \in \cP}\frac{1}{N}\sum_{i = 1}^N \ell(\boldsymbol{P}, \boldsymbol{\pi}_i).
\]

\paragraph{Object of study.} We aim to answer the standard generalization guarantee question: given a tolerance $\epsilon > 0$ and a failure probability $\delta \in (0, 1)$, what is the sample complexity $M(\epsilon, \delta)$ such that with probability at least $1 - \delta$ over the draw of problem instances $S = \{\boldsymbol{\pi}_1, \dots, \boldsymbol{\pi}_N\}$, where $N \geq M(\epsilon, \delta)$, we have 
$
    \bbE_{\boldsymbol{\pi} \sim \cD} [\ell(\hat{\boldsymbol{P}}_S, \boldsymbol{\pi})] \leq \bbE_{\boldsymbol{\pi} \sim \cD} [\ell(\boldsymbol{P}^*_\cD, \boldsymbol{\pi})] + \epsilon.
$ Consider the function class $\cL = \{\ell_{\boldsymbol{P}}: \Pi \rightarrow [-H, 0] \mid \boldsymbol{P} \in \cP\}$, where $\boldsymbol{\ell}_{\boldsymbol{P}}(\boldsymbol{\pi}) = \ell(\boldsymbol{P}, \boldsymbol{\pi})$. Theorem \ref{thm:pdim-to-generalization} suggests that the generalization guarantee is achievable by bounding the pseudo-dimension of $\cL$.

\section{Generalization Guarantee for Data-driven Input-Agnostic Projection Method for QPs} 
\label{sec:input-agnostic-QPs-guarantee}
This section provides the generalization guarantee for data-driven learning of the projection matrix $\boldsymbol{P}$ for QPs.

\subsection{Regularizing via Perturbing OQPs and the Perturbed Function Class}

There are two main obstacles to analyzing the generalization guarantee for data-driven learning the projection matrix in QPs. First, the optimal solution of QPs can lie arbitrarily anywhere in the feasible polyhedron. Second, when the matrix $\boldsymbol{Q}$ is singular, there can be infinitely many optimal solutions. To address this issue, we first introduce the perturbed function class $\cL_\gamma$, constructed by adding Tikhonov's regularization to the objective function of OQPs. This perturbation is equivalent to perturbing the matrix $\boldsymbol{Q}$ of the original input. After the perturbation: (i) the objective function of the perturbed OQPs and any perturbed PQPs becomes strictly convex, which favorably helps us localize the unique optimal solution and constructing the unrolled active set method; and (2) the perturbed function class $\cL_\gamma$ can approximate $\cL$ with arbitrary precision, and therefore analyzing the $\cL_\gamma$ can recover the guarantee for $\cL$. 

\begin{lemma} \label{lm:perturbed-QP-approximating-QP}
    Given a OQP $\boldsymbol{\pi} = (\boldsymbol{Q}, \boldsymbol{c}, \boldsymbol{A}, \boldsymbol{b})$, then there exists $\gamma > 0$ that is independent on $\boldsymbol{P}$, such that for the perturbed OQP $\boldsymbol{\pi}_\gamma = (\boldsymbol{Q}_\gamma, \boldsymbol{c}, \boldsymbol{A}, \boldsymbol{b})$, where $\boldsymbol{Q}_\gamma = \boldsymbol{Q} + \gamma \boldsymbol{I}_n$ and for any projection matrix $\boldsymbol{P} \in \cP$, we have 
    \[
    0 \leq \ell(\boldsymbol{P}, \boldsymbol{\pi}_\gamma) - \ell(\boldsymbol{P}, \boldsymbol{\pi}) \leq \frac{\gamma R^2}{2},\] 
    where $R$ comes from Assumption \ref{asmp:regularity}.3.
\end{lemma}
\proof
    Let $\boldsymbol{y}^*(\boldsymbol{P})$ is an optimal solution of the PQP, that is, $\ell(\boldsymbol{P}, \boldsymbol{\pi}) =  \frac{1}{2}{\boldsymbol{y}^*(\boldsymbol{P})}^\top \boldsymbol{P}^\top \boldsymbol{Q} \boldsymbol{P} \boldsymbol{y}^*(\boldsymbol{P}) + \boldsymbol{c}^\top \boldsymbol{P}\boldsymbol{y}^*(\boldsymbol{P})$, and let
    \begin{enumerate}
        \item[(i)] $f_{\boldsymbol{P}}(\boldsymbol{y}) = \frac{1}{2}\boldsymbol{y}^\top \boldsymbol{P}^\top \boldsymbol{Q} \boldsymbol{P} \boldsymbol{y} + \boldsymbol{c}^\top \boldsymbol{P}\boldsymbol{y}$ be the objective function of the PQP, and 
        \item[(ii)] $f_{\boldsymbol{P}, \boldsymbol{\gamma}}(\boldsymbol{y}) = \frac{1}{2}\boldsymbol{y}^\top \boldsymbol{P}^\top (\boldsymbol{Q} + \gamma\boldsymbol{I}_n)\boldsymbol{P} \boldsymbol{y} + \boldsymbol{c}^\top \boldsymbol{P}\boldsymbol{y} = f_{\boldsymbol{P}}(\boldsymbol{y}) + \frac{\gamma}{2}\|\boldsymbol{Py}\|_2^2$ be the objective function of the perturbed PQP. 
    \end{enumerate}
    Then $f_{\boldsymbol{P}}(\boldsymbol{y}^*(\boldsymbol{P})) = \ell(\boldsymbol{P}, \boldsymbol{\pi})$ by definition, and note that $\boldsymbol{y}^*(\boldsymbol{P})$ is a feasible point of the perturbed PQP $\boldsymbol{\pi}_\gamma$, meaning that $f_{\boldsymbol{P}, \gamma}(\boldsymbol{y}^*(\boldsymbol{P})) \geq \ell(\boldsymbol{P}, \boldsymbol{\pi}_\gamma)$. Besides, PQP and perturbed PQP have the same feasible region, with the objective of PQP $f_{\boldsymbol{P}}(\boldsymbol{y})$ is smaller than that of perturbed PQP $f_{\boldsymbol{P}, \gamma}(\boldsymbol{y})$, meaning that $\ell(\boldsymbol{P}, \boldsymbol{\pi}_\gamma) \geq \ell(\boldsymbol{P}, \boldsymbol{\pi})$. Combining the facts above, we have
    \begin{equation*}
        \begin{aligned}
            0 &\leq \ell(\boldsymbol{P}, \boldsymbol{\pi}_\gamma) - \ell(\boldsymbol{P}, \boldsymbol{\pi}) \leq f_{\boldsymbol{P}, \boldsymbol{\pi}}(\boldsymbol{y}^*(\boldsymbol{P})) - f_{\boldsymbol{P}}(\boldsymbol{y}^*(\boldsymbol{P})) \\
            &= \frac{\gamma}{2}\|\boldsymbol{P}\boldsymbol{y}^*(\boldsymbol{P})\|_2^2 \leq \frac{\gamma R^2}{2},
        \end{aligned}
    \end{equation*}
    where the final inequality comes from the fact that $\boldsymbol{P}\boldsymbol{y}^*(\boldsymbol{P})$ is a feasible point of OQP, and the feasible region is bounded by $R$ by Assumption \ref{asmp:regularity}.3. \qed 

    \begin{proposition} \label{prop:perturbed-OQP-favourable-structure}
        Any perturbed PQP of a perturbed OQP $\boldsymbol{\pi}_\gamma$ under a projection matrix $\boldsymbol{P}$ has a unique optimal solution. 
    \end{proposition}
    \proof
        First, notice that the matrix $\boldsymbol{P}^\top \boldsymbol{Q}_\gamma \boldsymbol{P}$ is positive definite.  To see that, for any $\boldsymbol{y} \in \bbR^k$ and $\boldsymbol{y} \neq \boldsymbol{0}$, we have $\boldsymbol{P}\boldsymbol{y} \neq \boldsymbol{0}_n$ since $\boldsymbol{P}$ is a full-column rank matrix. Therefore $\boldsymbol{y}^\top \boldsymbol{P}^\top \boldsymbol{Q}_\gamma \boldsymbol{P} \boldsymbol{y} = (\boldsymbol{P}\boldsymbol{y})^\top \boldsymbol{Q}_\gamma (\boldsymbol{P}\boldsymbol{y}) > 0$ as $\boldsymbol{Q}_\gamma = \boldsymbol{Q} + \gamma \boldsymbol{I}_n$ is a positive definite matrix. This implies that the objective value of the perturbed PQP $\boldsymbol{\pi}_\gamma$ is also strictly convex. Moreover, the perturbed PQP is feasible (admitting $\boldsymbol{0}_k$ as a feasible point by Proposition~\ref{prop:PQP-properties}) and bounded as $-H \leq \textup{OPT}(\boldsymbol{\pi}) \leq \ell(\boldsymbol{P}, \boldsymbol{\pi}) \leq \ell(\boldsymbol{P}, \boldsymbol{\pi}_\gamma)$. Therefore, the perturbed PQP has a unique optimal solution.
    \qed

    We now formally define the perturbed function class $\cL_\gamma$.
    \begin{definition}[Perturbed function class] \label{def:perturbed-function-class}
        Given $\gamma > 0$, the perturbed function class $\cL_\gamma$ is defined as  \[\cL_\gamma = \{\ell_{\boldsymbol{P}, \gamma}: \Pi \rightarrow [-H, 0] \mid \boldsymbol{P} \in \cP\},\] where $\ell_{\boldsymbol{P}, \gamma}(\boldsymbol{\pi}) = \ell(\boldsymbol{P}, \boldsymbol{\pi}_\gamma)$, $\boldsymbol{\pi}_\gamma = (\boldsymbol{Q} + \gamma \boldsymbol{I}_n, \boldsymbol{c}, \boldsymbol{A}, \boldsymbol{b})$, and $\ell(\boldsymbol{P}, \boldsymbol{\pi}_\gamma)$ is the optimal objective value of the perturbed PQP corresponding to the perturbed OQP $\boldsymbol{\pi}_\gamma$ and projection matrix $\boldsymbol{P}$.
    \end{definition}
    Next, we analyze the pseudo-dimension of $\cL_\gamma$, and later use the bound on the pseudo-dimension of $\cL_\gamma$ to bound the pseudo-dimension of $\cL$. Toward this goal, we need  to introduce several auxiliary results: In Section~\ref{sec:localizing-solution-perturbed-PQP}, we simplify the representation of the optimal solution using a localization technique, and in Section~\ref{sec:unrolled-active-set-method}, we devise a GJ algorithm to compute the optimal value of the perturbed PQP. These two results are the foundation for the bounding the psedo-dimension in Section~\ref{sec:recovering-guarantee-original-function-class}.

    \subsection{Localizing the  Solution of Perturbed PQPs}
    \label{sec:localizing-solution-perturbed-PQP}
      Consider the perturbed PQP corresponding to a projection matrix $\boldsymbol{P}$ and the perturbed OQP $\boldsymbol{\pi}_\gamma = (\boldsymbol{Q}_\gamma, \boldsymbol{c}, \boldsymbol{A}, \boldsymbol{b})$, and for convenience, let $\title{\boldsymbol{Q}} = \boldsymbol{P}^\top \boldsymbol{Q}_\gamma \boldsymbol{P}$, $\tilde{\boldsymbol{c}} = \boldsymbol{P}^\top \boldsymbol{c}$, and $\tilde{\boldsymbol{A}} = \boldsymbol{AP}$. We first prove that the solution of the perturbed PQP can be described using a simpler equality-constrained QP, of which the constraint matrix $\tilde{\boldsymbol{A}}_\cB$, extracted from the constraint matrix $\tilde{\boldsymbol{A}}$ of the perturbed PQP, has linearly independent rows. This serves as a localizing scheme for the optimal solution of a perturbed PQP, and is the foundation for the unrolled active set method (Algorithm \ref{alg:unrolled-active-set-method}) that we describe later.
    
    \begin{lemma} \label{lm:localizing-solution-pertubred-PQP}
        Let $\boldsymbol{y}^*$ be the (unique) optimal solution of perturbed PQP with the corresponding active set $\cA(\boldsymbol{y^*}) = \{i \in \{1, \dots, m\} \mid \tilde{\boldsymbol{A}}_i \boldsymbol{y}^* = \boldsymbol{b}_i\}$. Then there exists a subset $\cB \subseteq \cA(\boldsymbol{y}^*)$ such that:
        \begin{enumerate}
            \item The matrix $\tilde{\boldsymbol{A}}_\cB$ has linearly independent  row. Here $\boldsymbol{A}_\cB$ is the matrix formed by the row $i^{th}$ row of $\boldsymbol{A}$ for $i \in \cB$.
            \item $\boldsymbol{y}^*$ is the unique solution for the equality-constrained QP $\min_{\boldsymbol{y} \in \bbR^k} \left\{ \frac{1}{2}\boldsymbol{y}^\top\Tilde{\boldsymbol{Q}} \boldsymbol{y} + \Tilde{\boldsymbol{c}}^\top \boldsymbol{y} \mid \Tilde{\boldsymbol{A}}_\cB \boldsymbol{y} = \boldsymbol{b}_\cB\right\}$. 
        \end{enumerate}
    \end{lemma}

    \noindent \textit{Proof sketch.}
        The detailed proof can be found in Appendix \ref{apx:input-agnostic-QPs-guarantee}. Using KKT conditions, we claim that $-(\tilde{\boldsymbol{Q}}\boldsymbol{y}^* + \tilde{\boldsymbol{c}}) = \sum_{i \in \cA(\boldsymbol{y}^*)}\boldsymbol{\lambda}^*_i \cdot \tilde{\boldsymbol{A}}_i$. Note that $\sum_{i \in \cA(\boldsymbol{y}^*)}\boldsymbol{\lambda}^*_i \cdot \tilde{\boldsymbol{A}}_i$ is a conic combination, and using Conic's Carath\'{e}odory theorem (Proposition \ref{prop:conic-caratheodory}), we claim that there exists a subset $\cB \subseteq \cA(\boldsymbol{y}^*)$ such that there exists $\boldsymbol{\mu}_i \geq 0$ for $i \in \cB$ such that 
        $
            -(\tilde{\boldsymbol{Q}}\boldsymbol{y}^* + \tilde{\boldsymbol{c}}) = \sum_{j \in \cB}\boldsymbol{\mu}_j \cdot \boldsymbol{A}_j \Leftrightarrow  \tilde{\boldsymbol{Q}}\boldsymbol{y}^* + \tilde{\boldsymbol{c}} + \tilde{\boldsymbol{A}}_\cB^\top \boldsymbol{\mu}_\cB = \boldsymbol{0},
        $ 
        and $\tilde{\boldsymbol{A}}_{\cB}$ has linearly independent rows. Finally, we will show that $\boldsymbol{y}^*$ is the unique solution of the equality-constrained QP $\min_{\boldsymbol{y} \in \bbR^k} \left\{\frac{1}{2}\boldsymbol{y}^\top \tilde{\boldsymbol{Q}}\boldsymbol{y} + \tilde{\boldsymbol{c}}^\top \boldsymbol{y} \mid \tilde{\boldsymbol{A}}_\cB \boldsymbol{y} = \boldsymbol{b}_\cB\right\}$ by showing that $(\boldsymbol{y}^*, \boldsymbol{\mu}_\cB)$ is a KKT point of that problem.
    \qed

    
    \begin{remark} \label{rmk:localizing-solution-pertubred-PQP}
        In Lemma \ref{lm:localizing-solution-pertubred-PQP}, since $\tilde{\boldsymbol{A}}_\cB$ has linearly independent row and $\tilde{\boldsymbol{Q}}$ is positive definite, one can easily show that the KKT matrix $\boldsymbol{K} = \begin{bmatrix}
                    \tilde{\boldsymbol{Q}} & \tilde{\boldsymbol{A}}_\cB^\top \\
                    \tilde{\boldsymbol{A}}_\cB & \boldsymbol{0}
                \end{bmatrix}$  corresponding to the equality-constrained QP (Equation \ref{eq:equality-constrained-qp}) is invertible. Therefore, there exists $\boldsymbol{\lambda}_\cB$ such that $\begin{bmatrix}
                        \boldsymbol{y}^* \\
                        \boldsymbol{\lambda}^*_\cB
                    \end{bmatrix} = \boldsymbol{K}^{-1}\begin{bmatrix}
                        - \tilde{\boldsymbol{c}} \\
                        \boldsymbol{b}_\cB
                    \end{bmatrix}$. This point is helpful in designing the \textit{unrolled active set method} for computing $\ell(\boldsymbol{P}, \boldsymbol{\pi}_\gamma)$ in the next section.
    \end{remark}

    \subsection{The Unrolled Active Set Method} 
    \label{sec:unrolled-active-set-method}
    We now introduce the \textit{unrolled active set method}, a GJ algorithm with bounded complexities that exactly computes the optimal objective value of the perturbed PQP $\ell(\boldsymbol{P}, \boldsymbol{\pi}_\gamma)$ corresponding to the perturbed OQP $\boldsymbol{\pi}_\gamma$ and the projection matrix $\boldsymbol{P}$. 
      
    \noindent\textbf{Intuition.} The details of the unrolled active set method are demonstrated in Algorithm \ref{alg:unrolled-active-set-method}. Here, the algorithm is defined for each perturbed OQP $\boldsymbol{\pi}_\gamma$ and takes the projection matrix $\boldsymbol{P}$ as the input. The general idea is to check all the potential active subsets $\cA$ of rows of $\tilde{\boldsymbol{A}} = \boldsymbol{AP}$ up to $\min\{m, k\}$ elements. If we find $\cA$ such that KKT matrix $\boldsymbol{K} = \begin{bmatrix}
        \tilde{\boldsymbol{Q}} & \tilde{\boldsymbol{A}}_\cA^\top \\
        \tilde{\boldsymbol{A}}_\cA & \boldsymbol{0}.
    \end{bmatrix}$ is invertible, then we can use it to calculate the potential optimal solution $\boldsymbol{y}_\textup{cand}$ and Lagrangian $\boldsymbol{\lambda_\textup{cand}}$. We then check if $(\boldsymbol{y}_\textup{cand}, \boldsymbol{\lambda}_\textup{cand})$ is a KKT point of the perturbed PQP corresponding to the perturbed OQP $\boldsymbol{\pi}_\gamma$ and the projection matrix $\boldsymbol{P}$. If yes, then $\boldsymbol{y}_\textup{cand}$ is the optimal solution for the perturbed PQP, and we output the optimal objective value; else we move on to the next potential active subset $\cA$.

    \subsubsection{Correctness and GJ complexities.} In this section, we demonstrate that the unrolled active set method yields the optimal solution for the perturbed PQP. Then, we will show that the algorithm is also a GJ algorithm, and we will bound its predicate complexity and degree.
    
    \begin{lemma} \label{lm:unrolled-active-method-correctness}
        Given a perturbed OQP $\boldsymbol{\pi}_\gamma$, the algorithm $\Gamma_{\boldsymbol{\pi}_\gamma}$ described by Algorithm \ref{alg:unrolled-active-set-method} will output $\ell(\boldsymbol{P}, \boldsymbol{\pi}_\gamma)$.
    \end{lemma}
    \noindent \textit{Proof sketch.}
        The detailed proof can be found in Appendix~\ref{apx:input-agnostic-QPs-guarantee}. To prove the existence, showing that the algorithm guarantees to find an optimal solution $\boldsymbol{y}^*$, we have to use Lemma \ref{lm:localizing-solution-pertubred-PQP}, saying that there exists a subset $\cB \subset \cA(\boldsymbol{y}^*)$ such that $\boldsymbol{y}^*$ is the solution of the equality constrained QP corresponding to $\tilde{\boldsymbol{A}}_\cB$ with linearly independent rows. Then, we notice that the algorithm will check all subsets of $\{1, \dots, m\}$ of at most $\min(m, k)$ elements, hence it will eventually check $\cA = \cB$. When $\cA = \cB$, we verify that the candidate $\boldsymbol{y}_\textup{cand}$ and $\boldsymbol{\lambda}_\textup{cand}$ will pass all the primal and dual feasibility checks, and $\boldsymbol{y}_\textup{cand}$ is the optimal solution. For the correctness part, we will show that any $\boldsymbol{y}_\textup{cand}$ output by the algorithm is the optimal solution, by showing that $(\boldsymbol{y}_\textup{cand}, \boldsymbol{\lambda})$, where $\boldsymbol{\lambda}_\cA = \boldsymbol{\lambda}_\textup{cand}$ and $\boldsymbol{\lambda}_{\overline{\cA}} = \boldsymbol{0}$, is indeed a KKT point. 
    \qed 

    \begin{lemma} \label{lm:gj-complexity-unrolled-active-set} 
        Given a perturbed OQP $\boldsymbol{\pi}_\gamma$, the algorithm $\Gamma_{\boldsymbol{\pi}_\gamma}$ described by Algorithm \ref{alg:unrolled-active-set-method} is a GJ algorithm with degree $\cO(m + k)$ and predicate complexity $\cO\left(m\min\left(2^m, (\frac{em}{k})^k\right) \right)$.
    \end{lemma}
    \noindent \textit{Proof sketch.}
        The detailed proof is in Appendix \ref{apx:input-agnostic-QPs-guarantee}. First, note that $\tilde{\boldsymbol{Q}} = \boldsymbol{P}^\top \boldsymbol{Q}_\gamma \boldsymbol{P}$ is a matrix of which each entry is a polynomial in (the entries of) $\boldsymbol{P}$ of degree at most $2$. Similarly, each entry of $\tilde{\boldsymbol{A}} = \boldsymbol{AP}$ and $\tilde{\boldsymbol{c}} = \boldsymbol{P}^\top \boldsymbol{c}$ is a polynomial in $\boldsymbol{P}$ of degree at most $1$. We show that we have to check at most $\min(2^m, (em/k)^k)$ potential active sets. For each potential active set $\cA$, we show that the number of distinct predicates is $\cO(m)$ and the maximum degree of each predicate is $\cO(m + k)$. Combining those facts gives the result.
    \qed

    \begin{algorithm}[tb] 
        \caption{The unrolled active set method $\Gamma_{\boldsymbol{\pi}_\gamma}$ corresponding to the perturbed OQP $\boldsymbol{\pi}_\gamma = (\boldsymbol{Q}_\gamma, \boldsymbol{c}, \boldsymbol{A}, \boldsymbol{b})$}
            \label{alg:unrolled-active-set-method}
            \textbf{Input}: Projection matrix $\boldsymbol{P} \in \bbR^{n \times k}$\\
            \textbf{Output}: Optimal value of the perturbed PQP.
             \begin{algorithmic}[1]
            \STATE Set $\tilde{\boldsymbol{Q}} = \boldsymbol{P}^\top \boldsymbol{Q}_\gamma \boldsymbol{P}$, $\tilde{\boldsymbol{A}} = \boldsymbol{AP}$, and $\tilde{\boldsymbol{c}} = \boldsymbol{P}^\top \boldsymbol{c}$.
            
            \FOR{potential active set $\cA \subset \{1, \dots, m\}$, $\abs{\cA} \leq \min\{m, k\}$}
                \STATE Construct KKT matrix $\boldsymbol{K} = \begin{bmatrix}
                    \tilde{\boldsymbol{Q}} & \tilde{\boldsymbol{A}}_\cA^\top \\
                    \tilde{\boldsymbol{A}}_\cA & \boldsymbol{0}.
                \end{bmatrix}$
    
                \IF{$\det(\boldsymbol{K}) \neq 0$}
                    \STATE Compute $\begin{bmatrix}
                        \boldsymbol{y}_\textup{cand} \\
                        \boldsymbol{\lambda}_\textup{cand}
                    \end{bmatrix} = \boldsymbol{K}^{-1}\begin{bmatrix}
                        - \tilde{\boldsymbol{c}} \\
                        \boldsymbol{b}_\cA
                    \end{bmatrix}$.
    
                    \LineComment{Checking feasibility of potential solution $\boldsymbol{y}_\textup{cand}$}
                    \STATE $yFeasible = True$ 
                    \FOR{$j \not \in \cA$}
                        \IF{$\tilde{\boldsymbol{A}}_j^\top \boldsymbol{y}_\textup{cand} > \boldsymbol{b}_j$}
                            \STATE $yFeasible = False$
                            \STATE \textbf{break}
                        \ENDIF
                    \ENDFOR
                    \IF{$yFeasible$}
                        \LineComment{Checking validation of Lagrangian $\boldsymbol{\lambda}_{\textup{cand}}$}
                        $lambdaValid = True$
                        \FOR{$j \in \cA$}
                            \IF{$\boldsymbol{\lambda}_{\textup{cand}, j} < 0$}
                                \STATE $lambdaValid = False$
                                \STATE \textbf{break}
                            \ENDIF
                        \ENDFOR
                        
                        \IF{$lambdaValid$}
                            \STATE \textbf{return} $\frac{1}{2}\boldsymbol{y}_\textup{cand}^\top \tilde{\boldsymbol{Q}} \boldsymbol{y}_\textup{cand} + \tilde{\boldsymbol{c}}^\top \boldsymbol{y}_\textup{cand}$
                        \ENDIF
                    \ENDIF
                \ENDIF
            \ENDFOR
        
        \end{algorithmic}
    \end{algorithm}

\subsection{Pseudo-dimension Upper-bound Recovery for the Original Function Class} \label{sec:recovering-guarantee-original-function-class}
Using Lemma \ref{lm:gj-complexity-unrolled-active-set}, we now give a concrete upper-bound for the pseudo-dimension of the perturbed function class $\cL_\gamma$.
\begin{lemma}[Pseudo-dimension of $\cL_\gamma$]\label{lm:pdim-perturbed-bound}
     We have $\Pdim(\cL_\gamma) = \cO\left(nk \min(m, k \log m) \right)$ for any $\gamma > 0$.
\end{lemma}
Using the bound on the pseudo-dimension of perturbed function class $\Pdim(\cL_\gamma)$ and the connection between $\cL_\gamma$ and $\cL$ via Lemma \ref{lm:perturbed-QP-approximating-QP}, we can bound the the pseudo-dimension of the original function class $\cL$ as follows.
\begin{theorem}[Pseudo-dimension of $\cL$] \label{thm:pdim-original-bound}
    We have $\Pdim(\cL) = \cO\left(nk \min(m, k \log m) \right)$. 
\end{theorem}
\proof
    First, we claim that $0 \leq \textup{fatdim}_{\gamma R^2 / 2}\cL \leq \Pdim(\cL_\gamma)$ for any $\gamma > 0$, where $\textup{fatdim}_\alpha(\cL)$ is the fat-shattering dimension of $\cL$ at scale $\alpha$ (Definition \ref{def:fatdim}). To see that, assume $S = \{\boldsymbol{\pi}_1, \dots, \boldsymbol{\pi}_N\}$ is $\frac{\gamma R^2}{2}$ fat-shattered by $\cL$, meaning that there exists real-valued thresholds $r_1, \dots, r_N \in \bbR$ such that for any $I \subseteq \{1, \dots, N\}$, there exists $\ell_{\boldsymbol{P}} \in \cL$ such that
    \begin{align*}
        &f_{\boldsymbol{P}}(\boldsymbol{\pi}_i) > r_i + \frac{\gamma R^2}{2} \, \textup{for} \, i \in I, \quad \text{and} \\
    &f_{\boldsymbol{P}}(\boldsymbol{\pi}_j) < r_j - \frac{\gamma R^2}{2}\,\, \textup{for } j \not \in I.
    \end{align*}
    From Lemma \ref{lm:perturbed-QP-approximating-QP}, we have
    $0 \leq \ell_{\boldsymbol{P}, \gamma}(\boldsymbol{\pi}) - \ell_{\boldsymbol{P}}(\boldsymbol{\pi}) \leq \frac{\gamma R^2}{2}$ for any $\boldsymbol{\pi}$ and any $\boldsymbol{P} \in \cP$. This implies that $f_{\boldsymbol{P}, \gamma}(\boldsymbol{\pi}_i) > r_i$ if and only if $i \in I$. Therefore, $S$ is also pseudo-shattered by $\cL_{\gamma}$, which implies $0 \leq \textup{fatdim}_{\gamma R^2/ 2}(\cL) \leq \Pdim(\cL_\gamma)$. From Lemma~\ref{lm:pdim-perturbed-bound}, $\Pdim(\cL_\gamma) = \cO(nk\min(m, k\log m))$ for any $\gamma > 0$, therefore $0 \leq \textup{fatdim}_{\gamma R^2/2}(\cL) \leq C \cdot nk\min(m, k\log m)$ for any $\gamma > 0$ and some fixed constant $C$. Taking the limit $\gamma \rightarrow 0^+$ and using Proposition~\ref{prop:fatdim-properties}, we have $0 \leq \Pdim(\cL) \leq C\cdot nk\min(m, k\log m)$, or $\Pdim(\cL) = \cO(nk\min(m, k\log m)).$
\qed 

    Note that Theorem~\ref{thm:pdim-original-bound} is also applicable for data-driven learning projection matrix for LPs, as LP is a sub-problem of QP. Compared to the upper bound $\Pdim(\cL_\textup{LP}) = \cO(nk^2\log mk)$ by \citet[Theorem 4.4]{sakaue2024generalization}, our bound in Theorem \ref{thm:pdim-original-bound} is strictly tighter and applicable to both QPs and LPs.

\subsection{Lower Bound of Pseudo-dimension}
\label{sec:lower-bound}
For completeness, we also present the lower-bound for $\Pdim(\cL)$, of which the construction is inspired by the construction of learning projection matrix for LPs \cite{sakaue2024generalization}. See Appendix \ref{apx:input-agnostic-QPs-guarantee} for proof details.
\begin{proposition} \label{prop:lower-bound}
We have    $\Pdim(\cL) = \Omega(nk)$.
\end{proposition}

\section{Extension to Other Settings} \label{sec:extension}

This section explores several extensions of our results. 

\subsection{Learning to Match the Optimal Solution} \label{sec:match-solution}
In many cases, the decision-maker is interested in concrete solutions to implement. In these situations, the optimal value is less important than the optimal solution. We now consider the setting where we want to learn the projection matrix $\boldsymbol{P}$ such that the optimal solution of the PQP is close to that of the OQP. Such an approximative solution can also be used to warm-start an exact solver and accelerate the solving process. In this section, we propose an alternative objective value for learning $\boldsymbol{P}$. First, we assume the strict convexity of the problem instance, so that the optimal solution of the OQP is well-defined (unique).
\begin{assumption} \label{asmp:Q-PD}
    For any $\boldsymbol{\pi} = (\boldsymbol{Q}, \boldsymbol{c}, \boldsymbol{A}. \boldsymbol{b}) \in \Pi$, the matrix $\boldsymbol{Q}$ is positive definite.
\end{assumption}

Under such assumption, we seek the projection matrix $\boldsymbol{P}$ such that the recovered solution is close to the optimal solution of the PQP in expectation, i.e., 
\[
    \boldsymbol{P}^*_\cD \in \argmin_{\boldsymbol{P} \in \cP} \bbE_{\boldsymbol{\pi} \sim \cD} [\ell_\textup{match}(\boldsymbol{P}, \boldsymbol{\pi})],
\]
where 
$
    \ell_\textup{match}(\boldsymbol{P}, \boldsymbol{\pi}) = \|\boldsymbol{x}^*_{\boldsymbol{\pi}} - \boldsymbol{P}\boldsymbol{y}^*(\boldsymbol{P}, \boldsymbol{\pi})\|_2^2
$
is the matching loss, 
$\boldsymbol{x}^*_{\boldsymbol{\pi}} = \argmin_{\boldsymbol{x} \in \bbR^n}\left\{\frac{1}{2}\boldsymbol{x}^\top \boldsymbol{Q}\boldsymbol{x} + \boldsymbol{c}^\top \boldsymbol{x} \mid \boldsymbol{Ax} \leq \boldsymbol{b}\right\}$ is the optimal solution of the OQP, and $\boldsymbol{y}^*(\boldsymbol{P}, \boldsymbol{\pi}) = \argmin_{\boldsymbol{y} \in \bbR^k} \left\{\frac{1}{2}\boldsymbol{y}^\top\boldsymbol{P}^\top \boldsymbol{Q} \boldsymbol{P}\boldsymbol{y} + \boldsymbol{c}^\top \boldsymbol{P}\boldsymbol{y} \mid \boldsymbol{A\boldsymbol{P}}\boldsymbol{y} \leq \boldsymbol{b}\right\}$ is the optimal solution of the PQP. Again, since $\cD$ is unknown, we are instead given $N$ problem instances $S = \{\boldsymbol{\pi}_1, \dots, \boldsymbol{\pi}_N\}$ drawn i.i.d.~from $\cD$, and learn $\boldsymbol{P}$ via ERM
\[
    \hat{\boldsymbol{P}}_S \in \argmin_{\boldsymbol{P} \in \cP} \frac{1}{N} \sum_{i =1}^N \ell_\textup{match}(\boldsymbol{P}, \boldsymbol{\pi}_i).
\]
Let $\cL_\textup{match} = \{\ell_{\textup{match}, \boldsymbol{P}}: \Pi \rightarrow [-H, 0]\} \mid \boldsymbol{P} \in \cP \}$, where $\ell_{\textup{match}, \boldsymbol{P}}(\boldsymbol{\pi}) \coloneqq \ell_\textup{match}(\boldsymbol{P}, \boldsymbol{\pi})$. The following result provides the upper-bound for the pseudo-dimension of $\cL_\textup{match}$.
\begin{theorem} \label{thm:learning-to-match-solution}
    Assuming that all the QPs satisfies Assumption \ref{asmp:Q-PD} so that $\boldsymbol{x}^*_{\boldsymbol{\pi}}$ is defined uniquely. Then $\Pdim(\cL_\textup{match}) = \cO(nk \min(m, k \log m)).$
\end{theorem}
\noindent \textit{Proof sketch.} The detailed proof is presented in Appendix~\ref{apx:extension}. Given $\boldsymbol{\pi}$, the general idea is using a variant of the unrolled active set method (Algorithm \ref{alg:unrolled-active-set-method}) to calculate the optimal solution $\boldsymbol{y}^*(\boldsymbol{P}, \boldsymbol{\pi})$. Then $\ell_\textup{match}(\boldsymbol{P}, \boldsymbol{\pi})$ can also be calculated with a GJ algorithm with a bounded predicate complexity and degree, based on the GJ algorithm calculating $\boldsymbol{y}^*(\boldsymbol{P}, \boldsymbol{\pi})$. Finally, Theorem \ref{thm:gj} gives us the final guarantee.
\qed

\subsection{Input-aware Learning of Projection Matrix} \label{sec:input-aware}
In this section, we consider the setting of input-aware data-driven learning of the projection matrix for QPs, recently proposed by \citet{iwata2025learning} in the context of LP. Here, instead of learning a single projection matrix $\boldsymbol{P}$, we learn a mapping $f_{\boldsymbol{\theta}}: \Pi \rightarrow \cP$, e.g., a neural network, that takes a problem instance $\boldsymbol{\pi}$ drawn from $\cD$ and outputs the corresponding projection matrix $\boldsymbol{P}_{\boldsymbol{\pi}} = f_{\boldsymbol{\theta}}(\boldsymbol{\pi})$. With some computational trade-off for generating the projection matrix, this method has shown promising results, generating a better input-aware projection matrix that achieves better performance than an input-agnostic projection matrix while using the same projection dimension $k$. 

\noindent\textbf{Network architecture.} Inspired by \citet{iwata2025learning}, we assume that $f_{\boldsymbol{\theta}}$ is a neural network parameterized by $\boldsymbol{\theta} \in \Theta \subset \bbR^W$, where $W$ is the number of parameters of the neural network that takes the input $\boldsymbol{\pi}_\textup{flat}$ of size $n^2 + n + nm + m$ that is formed by flattening $\boldsymbol{Q}, \boldsymbol{c}, \boldsymbol{A}, \boldsymbol{b}$ in $\boldsymbol{\pi}$. Let $L$ be the number of hidden layers, and let $f_{\boldsymbol{\theta}}$ be the network of $L + 2$ layers, with the number of neurons of input layer is $W_0 = m^2 + n + mn + m$, that of the output layer is $W_{L + 2} = nk$, and that of $i^{th}$ layer is $W_i$ for $i \in \{1, \dots, L\}$. Each hidden layer uses ReLU as the non-linear activation function, and let $U = \sum_{i = 1}^NW_i$ be the number of hidden neurons. Consider the function class $\cL_{\textup{ia}} = \{\ell_{\boldsymbol{\theta}}: \Pi \rightarrow [-H, 0] \mid \boldsymbol{\theta} \in \Theta\}$, where
$
    \ell_{\boldsymbol{\theta}}(\boldsymbol{\pi}) \coloneqq \ell(f_{\boldsymbol{\theta}}(\boldsymbol{\pi}_\textup{flat}), \boldsymbol{\pi}).
$
Then we have the following result, which bounds the pseudo-dimension of $\cL_\textup{ia}$. The detailed proof is deferred to Appendix \ref{apx:extension}.
\begin{theorem} \label{thm:input-aware-setting}
    Assume that the output $f_{\boldsymbol{\theta}}(\boldsymbol{\pi})$ has full column rank, then $\Pdim(\cL_{\textup{ia}, \gamma}) = \cO(W(L\log(U + mk) + \min(m, k\log m)))$.
\end{theorem}

\section{Conclusion and Future Works}
We introduced the task of data-driven learning of a projection matrix for convex QPs. By a novel analysis approach, we establish the first upper bound on the pseudo-dimension of the learning projection matrix in QPs. Compared to the previous bound of \cite{sakaue2024generalization}, our new result is more general because it applies to both QPs and LPs and is strictly tighter. We further extend our analysis to learning to match the optimal solution and the input-aware setting. Our analysis opens many interesting directions for extension, including the conic programming and (mixed) integer programming.

\bibliography{references}

@inproceedings{bartlett2022generalization,
  title={Generalization bounds for data-driven numerical linear algebra},
  author={Bartlett, Peter and Indyk, Piotr and Wagner, Tal},
  booktitle={Conference on Learning Theory},
  pages={2013--2040},
  year={2022},
  organization={PMLR}
}

@article{pollard1984convergence,
  title={Convergence of Stochastic Processes},
  author={Pollard, David},
  journal={Springer Series in Statistics},
  year={1984},
  publisher={Springer New York}
}

@inproceedings{sakaue2023improved,
  title={Improved generalization bound and learning of sparsity patterns for data-driven low-rank approximation},
  author={Sakaue, Shinsaku and Oki, Taihei},
  booktitle={International Conference on Artificial Intelligence and Statistics},
  pages={1--10},
  year={2023},
  organization={PMLR}
}

@book{gass2003linear,
  title={Linear programming: methods and applications},
  author={Gass, Saul I},
  year={2003},
  publisher={Courier Corporation}
}

@article{sakaue2024generalization,
  title={Generalization bound and learning methods for data-driven projections in linear programming},
  author={Sakaue, Shinsaku and Oki, Taihei},
  journal={Advances in Neural Information Processing Systems},
  volume={37},
  pages={12825--12846},
  year={2024}
}

@inproceedings{iwata2025learning,
  title={Learning to generate projections for reducing dimensionality of heterogeneous linear programming problems},
  author={Iwata, Tomoharu and Sakaue, Shinsaku},
  booktitle={Forty-second International Conference on Machine Learning},
  year={2025}
}

@book{dostal2009optimal,
  title={Optimal Quadratic Programming Algorithms: With Applications to Variational Inequalities},
  author={Dost{\'a}l, Zdenek},
  volume={23},
  year={2009},
  publisher={Springer Science \& Business Media}
}

@article{amos2202tutorial,
  title={Tutorial on amortized optimization for learning to optimize over continuous domains},
  author={Amos, Brandon},
  journal={arXiv preprint arXiv:2202.00665},
  year={2022}
}

@article{huangfu2018parallelizing,
  title={Parallelizing the dual revised simplex method},
  author={Huangfu, Qi and Hall, JA Julian},
  journal={Mathematical Programming Computation},
  volume={10},
  number={1},
  pages={119--142},
  year={2018},
  publisher={Springer}
}

@article{chowdhury2022faster,
  title={Faster randomized interior point methods for tall/wide linear programs},
  author={Chowdhury, Agniva and Dexter, Gregory and London, Palma and Avron, Haim and Drineas, Petros},
  journal={Journal of Machine Learning Research},
  volume={23},
  number={336},
  pages={1--48},
  year={2022}
}

@article{applegate2021practical,
  title={Practical large-scale linear programming using primal-dual hybrid gradient},
  author={Applegate, David and D{\'\i}az, Mateo and Hinder, Oliver and Lu, Haihao and Lubin, Miles and O'Donoghue, Brendan and Schudy, Warren},
  journal={Advances in Neural Information Processing Systems},
  volume={34},
  pages={20243--20257},
  year={2021}
}

@article{d2020random,
  title={Random projections for quadratic programs},
  author={d’Ambrosio, Claudia and Liberti, Leo and Poirion, Pierre-Louis and Vu, Ky},
  journal={Mathematical Programming},
  volume={183},
  number={1},
  pages={619--647},
  year={2020},
  publisher={Springer}
}

@inproceedings{vu2019random,
  title={Random projections for quadratic programs over a Euclidean ball},
  author={Vu, Ky and Poirion, Pierre-Louis and d’Ambrosio, Claudia and Liberti, Leo},
  booktitle={International Conference on Integer Programming and Combinatorial Optimization},
  pages={442--452},
  year={2019},
  organization={Springer}
}

@article{vu2018random,
  title={Random projections for linear programming},
  author={Vu, Ky and Poirion, Pierre-Louis and Liberti, Leo},
  journal={Mathematics of Operations Research},
  volume={43},
  number={4},
  pages={1051--1071},
  year={2018},
  publisher={INFORMS}
}

@article{balcan2020data,
  title={Data-driven algorithm design},
  author={Balcan, Maria-Florina},
  journal={arXiv preprint arXiv:2011.07177},
  year={2020}
}

@article{gupta2020data,
  title={Data-driven algorithm design},
  author={Gupta, Rishi and Roughgarden, Tim},
  journal={Communications of the ACM},
  volume={63},
  number={6},
  pages={87--94},
  year={2020},
  publisher={ACM New York, NY, USA}
}

@article{indyk2019learning,
  title={Learning-based low-rank approximations},
  author={Indyk, Piotr and Vakilian, Ali and Yuan, Yang},
  journal={Advances in Neural Information Processing Systems},
  volume={32},
  year={2019}
}

@inproceedings{li2023learning,
  title={Learning the Positions in CountSketch},
  author={Li, Yi and Lin, Honghao and Liu, Simin and Vakilian, Ali and Woodruff, David},
  booktitle={The Eleventh International Conference on Learning Representations},
  year={2023}
}

@inproceedings{balcan2018learning,
  title={Learning to branch},
  author={Balcan, Maria-Florina and Dick, Travis and Sandholm, Tuomas and Vitercik, Ellen},
  booktitle={International conference on machine learning},
  pages={344--353},
  year={2018},
  organization={PMLR}
}

@inproceedings{goldberg1993bounding,
  title={Bounding the Vapnik-Chervonenkis dimension of concept classes parameterized by real numbers},
  author={Goldberg, Paul and Jerrum, Mark},
  booktitle={Proceedings of the Sixth Annual Conference on Computational Learning Theory},
  pages={361--369},
  year={1993}
}

@book{nocedal2006numerical,
  title={Numerical Optimization},
  author={Nocedal, Jorge and Wright, Stephen J},
  year={2006},
  publisher={Springer}
}

@inproceedings{dikin1967iterative,
  title={Iterative solution of problems of linear and quadratic programming},
  author={Dikin, Iliya Iosiphovich},
  booktitle={Doklady Akademii Nauk},
  volume={8},
  pages={674--675},
  year={1967}
}

@article{akchen2025column,
  title={Column-randomized linear programs: Performance guarantees and applications},
  author={Akchen, Yi-Chun and Misic, Velibor V},
  journal={Operations Research},
  volume={73},
  number={3},
  pages={1366--1383},
  year={2025},
  publisher={INFORMS}
}

@article{poirion2023random,
  title={Random projections of linear and semidefinite problems with linear inequalities},
  author={Poirion, Pierre-Louis and Lourenco, Bruno F and Takeda, Akiko},
  journal={Linear Algebra and its Applications},
  volume={664},
  pages={24--60},
  year={2023},
  publisher={Elsevier}
}

@article{chen2022learning,
  title={Learning to optimize: A primer and a benchmark},
  author={Chen, Tianlong and Chen, Xiaohan and Chen, Wuyang and Heaton, Howard and Liu, Jialin and Wang, Zhangyang and Yin, Wotao},
  journal={Journal of Machine Learning Research},
  volume={23},
  number={189},
  pages={1--59},
  year={2022}
}

@article{amos2023tutorial,
  title={Tutorial on amortized optimization},
  author={Amos, Brandon and others},
  journal={Foundations and Trends{\textregistered} in Machine Learning},
  volume={16},
  number={5},
  pages={592--732},
  year={2023},
  publisher={Now Publishers, Inc.}
}

@article{bengio2021machine,
  title={Machine learning for combinatorial optimization: A methodological tour d’horizon},
  author={Bengio, Yoshua and Lodi, Andrea and Prouvost, Antoine},
  journal={European Journal of Operational Research},
  volume={290},
  number={2},
  pages={405--421},
  year={2021},
  publisher={Elsevier}
}

@article{sauer1972density,
  title={On the density of families of sets},
  author={Sauer, Norbert},
  journal={Journal of Combinatorial Theory, Series A},
  volume={13},
  number={1},
  pages={145--147},
  year={1972},
  publisher={Elsevier}
}

@book{horn2012matrix,
  title={Matrix Analysis},
  author={Horn, Roger A and Johnson, Charles R},
  year={2012},
  publisher={Cambridge University Press}
}

@inproceedings{bartlett1994fat,
  title={Fat-shattering and the learnability of real-valued functions},
  author={Bartlett, Peter L and Long, Philip M and Williamson, Robert C},
  booktitle={Proceedings of the Seventh Annual Conference on Computational Learning Theory},
  pages={299--310},
  year={1994}
}

@article{cheng2024sample,
  title={Sample complexity of algorithm selection using neural networks and its applications to branch-and-cut},
  author={Cheng, Hongyu and Khalife, Sammy and Fiedorowicz, Barbara and Basu, Amitabh},
  journal={Advances in Neural Information Processing Systems},
  volume={37},
  pages={25036--25060},
  year={2024}
}

@book{anthony2009neural,
  title={Neural Network Learning: Theoretical Foundations},
  author={Anthony, Martin and Bartlett, Peter L},
  year={2009},
  publisher={Cambridge University Press}
}

@article{warren1968lower,
  title={Lower bounds for approximation by nonlinear manifolds},
  author={Warren, Hugh E},
  journal={Transactions of the American Mathematical Society},
  volume={133},
  number={1},
  pages={167--178},
  year={1968},
  publisher={JSTOR}
}

@article{milgrom2002envelope,
  title={Envelope theorems for arbitrary choice sets},
  author={Milgrom, Paul and Segal, Ilya},
  journal={Econometrica},
  volume={70},
  number={2},
  pages={583--601},
  year={2002},
  publisher={Wiley Online Library}
}

@article{balcan2022provably,
  title={Provably tuning the ElasticNet across instances},
  author={Balcan, Maria-Florina F and Khodak, Misha and Sharma, Dravyansh and Talwalkar, Ameet},
  journal={Advances in Neural Information Processing Systems},
  volume={35},
  pages={27769--27782},
  year={2022}
}

@article{balcan2023new,
  title={New bounds for hyperparameter tuning of regression problems across instances},
  author={Balcan, Maria-Florina F and Nguyen, Anh and Sharma, Dravyansh},
  journal={Advances in Neural Information Processing Systems},
  volume={36},
  pages={80066--80078},
  year={2023}
}

@article{balcan2025algorithm,
  title={Algorithm configuration for structured Pfaffian settings},
  author={Balcan, Maria Florina and Nguyen, Anh Tuan and Sharma, Dravyansh},
  journal={Transactions on Machine Learning Research}, 
  year={2025}
}

@inproceedings{
    balcan2025sample,
    title={Sample complexity of data-driven tuning of model hyperparameters in neural networks with structured parameter-dependent dual function},
    author={Maria Florina Balcan and Anh Tuan Nguyen and Dravyansh Sharma},
    booktitle={The Thirty-ninth Annual Conference on Neural Information Processing Systems},
    year={2025},
    url={https://openreview.net/forum?id=pFFFRi2TcC}
}

\appendix

\section{Additional Backgrounds on Learning Theory}

In this section, we will go through the definition of the fat-shattering dimension of a real-valued function class and its connection to the pseudo-dimension. This definition of learning-theoretic complexity is useful in our case, when we want to draw the connection between the pseudo-dimension of the perturbed function class $\cL_\gamma$ and the pseudo-dimension of the original function class $\cL$, as in Section \ref{sec:recovering-guarantee-original-function-class}.

\begin{definition}[Fat-shattering dimension, \cite{bartlett1994fat}] \label{def:fatdim}
     Consider a real-valued function class $\cL$, of which each function $\ell$ takes input $\boldsymbol{\pi}$ in $\Pi$ and output $\ell(\boldsymbol{\pi}) \in [-H, 0]$. Given a set of inputs $S = \{\boldsymbol{\pi_1}, \dots, \boldsymbol{\pi}_N\} \subset \Pi$, we say that $S$ is fat-shattered at scale $\alpha > 0$ if there exists real-valued thresholds $r_1, \dots, r_N \in \bbR$ such that for any index $I \subseteq \{1, \dots, M\}$, there exists $\ell \in \cL$ such that
     \[
        f(\boldsymbol{\pi}_i) > r_i + \alpha \, \text{ for } i\in I, \text{ and } f(\boldsymbol{\pi}_j) < r_j - \alpha \,\text{ for }  \, j \not \in I.
     \]
    The fat-shattering dimension of $\cL$ at scale $\alpha$, denote $\textup{fatdim}_\alpha(\cL)$ is the the size of the largest set $S$ that can be shattered at scale $\alpha$ by $\cL$.
\end{definition}

The following results demonstrate some basic property of fat-shattering dimension and its connection to the pseudo-dimension. 
\begin{proposition}[\cite{bartlett1994fat}] \label{prop:fatdim-properties}
    Let $\cL$ be a real-valued function class, then:
    \begin{enumerate}
        \item For all $\alpha > 0$, $\textup{fatdim}_\alpha(\cL) \leq \Pdim(\cL)$.
        \item The function $\textup{fatdim}_\alpha(\cL)$ is non-decreasing with $\alpha$.
        \item If a finite set $S$ is pseudo-shattered, then there is some $\alpha_0 > 0$ such that for all $\alpha < \alpha_0$, the set $S$ is fat-shattered at scale $\alpha$.
        \item $\lim_{\alpha \rightarrow 0^+} \textup{fatdim}_\alpha(\cL) = \Pdim(\cL)$.
    \end{enumerate}
    
\end{proposition}

\section{Proofs for Section \ref{sec:input-agnostic-QPs-guarantee}} \label{apx:input-agnostic-QPs-guarantee}

We first recall the Sauer-Shelah lemma, which is a well-known result in combinatorics that allows us to bound the sum of a combinatorial sequence. 
\begin{lemma}[Sauer-Shelah lemma \cite{sauer1972density}] \label{lm:sauer-shelah}
     Let $1 \leq k \leq n$, where $k$ and $n$ are positive integers. Then $$\sum_{j=0}^k{n \choose j} \leq \left(\frac{en}{k}\right)^k.$$
\end{lemma}

We then recall the Warren's theorem \cite{warren1968lower}, which bounds the number of sign patterns that a sequence of polynomials with bounded degrees can create.
\begin{lemma}[Warren's theorem \cite{warren1968lower}] \label{lm:warren}
    Let $p_1(\boldsymbol{x}), \dots, p_m(\boldsymbol{x})$ be polynomials in $n$ variables of degree at most $d$. Then the number of sign patterns
    \[
        (\sign(p_1(\boldsymbol{x})), \dots, \sign(p_m(\boldsymbol{x})))
    \]
    acquired by varying $\boldsymbol{x}$ is at most $\left(\frac{8edm}{n}\right)^n$
    
\end{lemma}

    \subsection{Proofs for Section \ref{sec:localizing-solution-perturbed-PQP}}
    \label{apx:localizing-solution-perturbed-PQP}

    We now present the formal proof for Lemma \ref{lm:localizing-solution-pertubred-PQP}.
    

    \proof[Proof of Lemma \ref{lm:localizing-solution-pertubred-PQP}]
        Since $\boldsymbol{y}^*$ is the optimal solution of the perturbed PQP problem, then there exists a vector $\boldsymbol{\lambda}^* \in \bbR^m$ such that $(\boldsymbol{y}^*, \boldsymbol{\lambda}^*)$ that satisfies the KKT conditions:
        \begin{enumerate}
            \item Stationarity: $\tilde{\boldsymbol{Q}}\boldsymbol{y}^* + \tilde{\boldsymbol{c}} + \tilde{\boldsymbol{A}}^\top \boldsymbol{\lambda}^* = \boldsymbol{0}$.
            \item Primal feasibility: $\tilde{\boldsymbol{A}}\boldsymbol{y}^* \leq \boldsymbol{b}$. \item Dual feasibility: $\boldsymbol{\lambda}^* \geq 0$.
            \item Complementary slackness: $\boldsymbol{\lambda}^*_i(\tilde{\boldsymbol{A}}_i \boldsymbol{y}^* - \boldsymbol{b}_i) = 0$, for $i \in \{1, \dots, m\}$.
        \end{enumerate}
        From the property of the active set $\cA(\boldsymbol{y}^*)$ and the complementary slackness property, we have $\boldsymbol{\lambda}^*_{\overline{\cA}(\boldsymbol{y}^*)} = \boldsymbol{0}$, where $\overline{\cA}(\boldsymbol{y}^*) = \{1, \dots, m\} \setminus \cA(\boldsymbol{y}^*)$ is the complement of the active set $\cA(\boldsymbol{y}^*)$. Combining the fact that $\boldsymbol{\lambda}^*_{\overline{\cA}(\boldsymbol{y}^*)} = \boldsymbol{0}$ and the stationary condition above, we have
        \begin{equation*}
            \begin{aligned}
                &\tilde{\boldsymbol{Q}}\boldsymbol{y}^* + \tilde{\boldsymbol{c}} + \tilde{\boldsymbol{A}}^\top_{\cA(\boldsymbol{x}^*)} \boldsymbol{\lambda}^*_{\cA(\boldsymbol{y}^*)} = \boldsymbol{0} \\
                \Rightarrow &-(\tilde{\boldsymbol{Q}}\boldsymbol{y}^* + \tilde{\boldsymbol{c}}) = \sum_{i \in \cA(\boldsymbol{y}^*)}\boldsymbol{\lambda}^*_i \cdot \tilde{\boldsymbol{A}}_i,
            \end{aligned}
        \end{equation*}
        where $\tilde{\boldsymbol{A}}_i$ is the $i^{th}$ row of $\tilde{\boldsymbol{A}}$. Since $\boldsymbol{\lambda}^*_i \geq 0$ for all $i \in \{1, \dots, m\}$, we can see that $-(\tilde{\boldsymbol{Q}}\boldsymbol{y}^* + \tilde{\boldsymbol{c}})$ is the conic combination of $\tilde{\boldsymbol{A}}_i$ for $i \in \cA(\boldsymbol{y}^*)$. We now recall the Conic's Carath\'{e}odory theorem, which can simplify the representation of a conic combination.

        \begin{proposition}[Conic's Carath\'{e}odory theorem] \label{prop:conic-caratheodory}
            If $\boldsymbol{v} \in \bbR^n$ lies in $\textup{Conic}(\boldsymbol{S})$, where $\boldsymbol{S} = \{\boldsymbol{s}_1, \dots, \boldsymbol{s}_t\} \subset \bbR^n$, then $\boldsymbol{v}$ can be rewritten as a linear combination of at most $n$ linearly independent vector from $\boldsymbol{S}$.
        \end{proposition}

        Using  the Conic's Carath\'{e}odory theorem, we claim that there exists a index set $\cB \subset \cA(\boldsymbol{y}^*)$, and $\boldsymbol{\mu}_j \geq 0$ for $j \in \cB$ such that  
        \begin{equation} \label{eq:conic-form}
            -(\tilde{\boldsymbol{Q}}\boldsymbol{y}^* + \tilde{\boldsymbol{c}}) = \sum_{j \in \cB}\boldsymbol{\mu}_j \cdot \boldsymbol{A}_j \Leftrightarrow  \tilde{\boldsymbol{Q}}\boldsymbol{y}^* + \tilde{\boldsymbol{c}} + \tilde{\boldsymbol{A}}_\cB^\top \boldsymbol{\mu}_\cB = \boldsymbol{0}.
        \end{equation}
        Now, consider the new equality-constrained QP 
        \begin{equation} \label{eq:equality-constrained-qp}
            \min_{\boldsymbol{y} \in \bbR^k} \left\{\frac{1}{2}\boldsymbol{y}^\top \tilde{\boldsymbol{Q}}\boldsymbol{y} + \tilde{\boldsymbol{c}}^\top \boldsymbol{y} \mid \tilde{\boldsymbol{A}}_\cB \boldsymbol{y} = \boldsymbol{b}_\cB\right\},
        \end{equation}
        and we claim that $\boldsymbol{y}^*$ is the (unique) solution of the problem above, by claiming that $(\boldsymbol{y}^*, \boldsymbol{\mu}_\cB)$ is a KKT point of the equality-constrained QP. 
        \begin{enumerate}
            \item First, from Equation~\eqref{eq:conic-form}, we have $\tilde{\boldsymbol{Q}}\boldsymbol{y}^* + \tilde{\boldsymbol{c}} + \tilde{\boldsymbol{A}}_\cB^\top \boldsymbol{\mu}_\cB = \boldsymbol{0}$. Therefore $(\boldsymbol{y}^*, \boldsymbol{\mu}_\cB)$ satisfies the stationarity condition.
            \item Since $\cB \subseteq \cA(\boldsymbol{y}^*)$, then $\tilde{\boldsymbol{A}}_i^\top \boldsymbol{y}^* = \boldsymbol{b}_i$ for $i \in \cB$. Therefore $\boldsymbol{y}^*$ satisfies the primal feasibility constraints. 
            \item The dual feasibility and complementary slackness are satisfied since this is an equality-constrained QP. 
        \end{enumerate}
        Therefore, $(\boldsymbol{y}^*, \boldsymbol{\mu}_\cB)$ is a KKT point of the equality-constrained QP and therefore $\boldsymbol{y}^*$ is an optimal solution. Moreover, since the objective function of the equality-constrained QP is strictly convex, $\boldsymbol{y}^*$ is the unique optimal solution. 
    \qed 
    \subsection{Proofs for Section \ref{sec:unrolled-active-set-method}} 
    \label{apx:unrolled-active-set-method}

    We now present the formal proof of Lemma \ref{lm:unrolled-active-method-correctness}, which shows the correctness for the unrolled active set method in Algorithm~\ref{alg:unrolled-active-set-method}.

    \proof[Proof of Lemma~\ref{lm:unrolled-active-method-correctness}]
        \textbf{Existence.} We will first show that $\Gamma_{\boldsymbol{\pi}_\gamma}$ guarantees to find an optimal solution $\boldsymbol{y}^*$ for the perturbed PQP corresponding to the perturbed OQP $\boldsymbol{\pi}_\gamma$ and the input projection matrix $\boldsymbol{P}$. From Lemma \ref{lm:localizing-solution-pertubred-PQP}, there exists $\cB \subseteq \cA(\boldsymbol{y}^*)$ such that $\tilde{\boldsymbol{A}}_\cB$ has linearly dependent rows, and $\boldsymbol{y}^*$ is the solution of the equality-constrained problem
        \[
            \min_{\boldsymbol{y} \in \bbR^k} \left\{\frac{1}{2}\boldsymbol{y}^\top \tilde{\boldsymbol{Q}}\boldsymbol{y} + \tilde{\boldsymbol{c}}^\top \boldsymbol{y} \mid \tilde{\boldsymbol{A}}_\cB \boldsymbol{y} = \boldsymbol{b}_\cB\right\},
        \]
        where $\tilde{\boldsymbol{Q}} = \boldsymbol{P}^\top \boldsymbol{Q}_\gamma \boldsymbol{P}$, $\tilde{\boldsymbol{c}} = \boldsymbol{P}^\top \boldsymbol{c}$, and $\tilde{\boldsymbol{A}} = \boldsymbol{AP}$. Since Algorithm \ref{alg:unrolled-active-set-method} will check all $\cA \subset \{1, \dots, m\}$ and $\abs{\cA} \leq k$, the algorithm $\Gamma_{\boldsymbol{\pi}_\gamma}$ will eventually select $\cA = \cB$. When $\Gamma_{\boldsymbol{\pi}_\gamma}$ selects $\cA = \cB$:
        \begin{enumerate}
            \item The KKT matrix $\boldsymbol{K} = \begin{bmatrix}
                \tilde{\boldsymbol{Q}} & \tilde{\boldsymbol{A}}_\cA^\top \\
                \tilde{\boldsymbol{A}}_\cA & \boldsymbol{0}.
            \end{bmatrix}$ is invertible, since $\tilde{\boldsymbol{A}}_\cA$ has linearly independent rows, and $\tilde{\boldsymbol{Q}}$ is positive definite.
            \item Then $\Gamma_{\boldsymbol{\pi}_\gamma}$ computes $\begin{bmatrix}
                        \boldsymbol{y}_\textup{cand} \\
                        \boldsymbol{\lambda}_\textup{cand}
                    \end{bmatrix} = \boldsymbol{K}^{-1}\begin{bmatrix}
                        - \tilde{\boldsymbol{c}} \\
                        \boldsymbol{b}_\cA
                    \end{bmatrix}$. From Lemma \ref{lm:localizing-solution-pertubred-PQP} and Remark \ref{rmk:localizing-solution-pertubred-PQP}, $\boldsymbol{y}_\textup{cand}$ is the optimal solution of the perturbed PQP corresponding to $\pi_\gamma$ and $\boldsymbol{P}$.
            \item Since $\boldsymbol{y}_\textup{cand}$ is the optimal solution, then the KKT conditions check will automatically pass.
        \end{enumerate}
        Therefore, $\boldsymbol{y}_\textup{cand}$ is the optimal solution, and $\Gamma_{\boldsymbol{\pi}_\gamma}$ will return the optimal value $\ell(\boldsymbol{P}, \boldsymbol{\pi}_\gamma)$ for the perturbed PQP.

        \textbf{Correctness.} We then show that any value $\boldsymbol{y}_\textup{cand}$ that $\Gamma_{\boldsymbol{\pi}_\gamma}$ (with the corresponding value $\frac{1}{2}\boldsymbol{y}_\textup{cand}^\top \tilde{\boldsymbol{Q}} \boldsymbol{y}_\textup{cand} + \tilde{\boldsymbol{c}}^\top \boldsymbol{y}_\textup{cand}$ is indeed the optimal solution for the perturbed PQP. To do that, we just have to verify $(\boldsymbol{y}_\textup{cand}, \boldsymbol{\lambda})$, where $\boldsymbol{\lambda}_\cA = \boldsymbol{\lambda}_\textup{cand}$ calculated by the algorithm, and $\boldsymbol{\lambda}_{\overline{\cA}} = \boldsymbol{0}$, satisfies the KKT conditions of the perturbed PQP. Here $\cA$ is the potential active set corresponding to $\boldsymbol{y}_\textup{cand}$ and $\overline{\cA} = \{1, \dots, m\} \setminus \cA$. 
        \begin{enumerate}
            \item From Algorithm \ref{alg:unrolled-active-set-method}, $\begin{bmatrix}
                        \boldsymbol{y}_\textup{cand} \\
                        \boldsymbol{\lambda}_\textup{cand}
                    \end{bmatrix} = \boldsymbol{K}^{-1}\begin{bmatrix}
                        - \tilde{\boldsymbol{c}} \\
                        \boldsymbol{b}_\cA
                    \end{bmatrix}$, meaning that $\tilde{\boldsymbol{Q}}\boldsymbol{y}_\textup{cand} + \tilde{\boldsymbol{c}} + \tilde{\boldsymbol{A}}_\cA^\top \boldsymbol{\lambda}_\textup{cand} = \boldsymbol{0}.$ And note that $\boldsymbol{\lambda}_\cA = \boldsymbol{\lambda}_\cA$ and $\boldsymbol{\lambda}_{\overline{\cA}} = \boldsymbol{0}$ by the definition above, we have $\tilde{\boldsymbol{Q}}\boldsymbol{y}_\textup{cand} + \tilde{\boldsymbol{c}} + \tilde{\boldsymbol{A}}^\top \boldsymbol{\lambda} = \boldsymbol{0}$, meaning that $(\boldsymbol{y}_\textup{cand}, \boldsymbol{\lambda})$ satisfies the stationarity condition.

            \item From Algorithm \ref{alg:unrolled-active-set-method}, $\boldsymbol{y}_\textup{cand}$ passes the feasibility check, meaning that it satisfies the primal feasibility condition.
            \item From Algorithm \ref{alg:unrolled-active-set-method}, $\boldsymbol{\lambda}_\textup{cand, i} \geq 0$ for all $i \in \cA$, and by definition $\boldsymbol{\lambda}_j = 0$ for all $j \in \overline{\cA}$. Therefore $\boldsymbol{\lambda}$ satisfies the dual feasibility condition.
            \item For $i \in \cA$, we have $\boldsymbol{\lambda}_i \cdot (\tilde{\boldsymbol{A}}_i^\top \boldsymbol{y}_\textup{cand} - \boldsymbol{b}_i) = 0$ since $\tilde{\boldsymbol{A}}_i^\top \boldsymbol{y}_\textup{cand} - \boldsymbol{b}_i$ from the property of active set. For $i \in \overline{\cA}$, $\boldsymbol{\lambda}_i \cdot (\tilde{\boldsymbol{A}}_i^\top \boldsymbol{y}_\textup{cand} - \boldsymbol{b}_i) = 0$ since $\boldsymbol{\lambda}_i = 0$ by definition. Therefore $(\boldsymbol{y}_\textup{cand}, \boldsymbol{\lambda})$ satisfies the complementary slackness.
        \end{enumerate}
        Therefore, $(\boldsymbol{y}_\textup{cand}, \boldsymbol{\lambda})$ is indeed a KKT point of the perturbed PQP, therefore $\boldsymbol{y}_\textup{cand}$ is its optimal solution and $\ell(\boldsymbol{P}, \boldsymbol{\pi}_\gamma) = \frac{1}{2}\boldsymbol{y}_\textup{cand}^\top \tilde{\boldsymbol{Q}} \boldsymbol{y}_\textup{cand} + \tilde{\boldsymbol{c}}^\top \boldsymbol{y}_\textup{cand}$.
    \qed 

    We now present the formal proof of Lemma \ref{lm:gj-complexity-unrolled-active-set}, which shows that the unrolled active set method is a GJ algorithm, and thereby bounds the predicate complexity and degree of the algorithm.

\proof[Proof of Lemma~\ref{lm:gj-complexity-unrolled-active-set}]
        First, note that $\tilde{\boldsymbol{Q}} = \boldsymbol{P}^\top \boldsymbol{Q}_\gamma \boldsymbol{P}$ is a matrix of which each entry is a polynomial in (the entries of) $\boldsymbol{P}$ of degree at most $2$. Similarly, each entry of $\tilde{\boldsymbol{A}} = \boldsymbol{AP}$ and $\tilde{\boldsymbol{c}} = \boldsymbol{P}^\top \boldsymbol{c}$ is a polynomial in $\boldsymbol{P}$ of degree at most $1$.

        Let $t = \min \{m, k \}$. From the algorithm, we have to consider (in the worst case) all subsets $\cA$ of $\{1, \dots, m\}$ with at most $t$ elements. Therefore, we have to consider at most $\min\left\{2^m, \left(\frac{em}{k}\right)^k \right\}$ subsets, where $2^m$ corresponds to the case $m \leq k$ and $\left(\frac{em}{k}\right)^k$ corresponds to the case $k < m$ and using Sauer-Shelah lemma (Lemma \ref{lm:sauer-shelah}).

        For each potential active set $\cA$: 
        \begin{enumerate}
            \item We have to check if $\det(\boldsymbol{K}) \neq 0$. Since $\boldsymbol{K} = \begin{bmatrix}
                    \tilde{\boldsymbol{Q}} & \tilde{\boldsymbol{A}}_\cA^\top \\
                    \tilde{\boldsymbol{A}}_\cA & \boldsymbol{0}.
                \end{bmatrix}$, each entry of $\boldsymbol{K}$ is a polynomial in $\boldsymbol{P}$ of degree at most $2$, and the size of $\boldsymbol{K}$ is $(k + \abs{\cA}) \times (k + \abs{\cA})$. Therefore, $\det(\boldsymbol{K})$ is a polynomial in $\boldsymbol{P}$ of degree at most $2(k + t)$. Besides, we have to check $\det(\boldsymbol{K}) \neq 0$ by checking $\det(\boldsymbol{K}) \geq 0$ and $-\det(\boldsymbol{K}) \geq 0$, which creates two distinct predicates.

            \item If $\det(\boldsymbol{K}) \neq 0$, we calculate $\boldsymbol{K}^{-1}$ via adjugate matrix, i.e., $\boldsymbol{K}^{-1} = \frac{\textup{adj}(\boldsymbol{K})}{\det(\boldsymbol{K})}$ \cite{horn2012matrix}. Therefore, each entry of $\boldsymbol{K}^{-1}$ is a rational function of $\boldsymbol{P}$ of degree at most $2(k + t)$. Then note that $\begin{bmatrix}
                        \boldsymbol{y}_\textup{cand} \\
                        \boldsymbol{\lambda}_\textup{cand}
                    \end{bmatrix} = \boldsymbol{K}^{-1}\begin{bmatrix}
                        - \tilde{\boldsymbol{c}} \\
                        \boldsymbol{b}_\cA
                    \end{bmatrix}$, meaning that each entry of $\boldsymbol{y}_\textup{cand}$ and $\boldsymbol{\lambda}_\textup{cand}$ is a rational function of $\boldsymbol{P}$ of degree at most $2(k + t) + 1$.

            \item After acquiring $(\boldsymbol{y}_\textup{cand}, \boldsymbol{\lambda}_\textup{cand})$, we have to check the primal and dual feasibility, which requires $m - \abs{\cA}$ distinct predicates of degree at most $2(k + t) + 2$ for checking $\tilde{\boldsymbol{A}}_j^\top \boldsymbol{y}_\textup{cand} \leq \boldsymbol{b}_j$, and $\abs{\cA}$ distinct predicates of degree at most $2(k + t) + 1$ for checking $\boldsymbol{\lambda}_{\textup{cand}, j} \geq 0$. The total distinct predicates in each steps is $\cO(m)$
        \end{enumerate}
        In total, in every steps, $\Gamma_{\boldsymbol{\pi}_\gamma}$ invovles in rational functions of $\boldsymbol{P}$, hence it is a GJ algorithm. For each active set $\cA$, the algorithm creates $\cO(m)$ distinct predicates of degree at most $\cO(k + t) \leq \cO(m + k)$. Therefore, the degree of $\boldsymbol{\Gamma}_{\boldsymbol{\pi}_\gamma}$ is $\Delta = \cO(m + k)$ and the predicate complexity is $\Lambda = \cO\left(m\min\left(2^m, (\frac{em}{k})^k\right) \right)$.
    \qed

\subsection{Proofs for Section \ref{sec:lower-bound}}
In this section, we will provide the detailed construction for the lower-bound presented in Proposition~\ref{prop:lower-bound}. The idea of the construction is already presented in the work by \cite{sakaue2024generalization} in the context of LPs. We adapt this approach to the case of QPs.

\proof[Proof of Proposition~\ref{prop:lower-bound}]
    We will construct the lower-bound by constructing a set of $(n - 2k)k$ QP instances that $\cL$ can shatter. For $r = 1, \dots, n - 2k$ and $s = 1, \dots, k$, we consider QP problem instance $\boldsymbol{\pi}_{r,s} = (\boldsymbol{Q}, \boldsymbol{c}_r, \boldsymbol{A}, \boldsymbol{b}_s) \in \Pi$, where
    \[
        \boldsymbol{Q} = \boldsymbol{0}_{n \times n}, \boldsymbol{c}_r = \begin{bmatrix}
            \mathbf{e}_{r} \\
            \boldsymbol{0}_{2k}
        \end{bmatrix}, \boldsymbol{A} = \begin{bmatrix}
            \boldsymbol{0}_{2k, n - 2k} & \boldsymbol{I}_{2k}
        \end{bmatrix},
        \boldsymbol{b}_s = \begin{bmatrix}
            \mathbf{e}_s \\
            \boldsymbol{0}_k
        \end{bmatrix},
    \]
    and $\mathbf{e}_r$ and $\mathbf{e}_s$ are the $r^{th}$ and $s^{th}$ standard basis vectors of $\bbR^{n - 2k}$ and $\bbR^{k}$, respectively. We consider the functions $\ell_{\boldsymbol{P}}: \Pi \rightarrow \bbR$, where the projection matrix $\boldsymbol{P}$ takes the form
    \[
        \boldsymbol{P} = \begin{bmatrix}
            \boldsymbol{T} \\
            \boldsymbol{I}_k \\
            -\boldsymbol{I}_{k},
        \end{bmatrix}
    \]
    and $\boldsymbol{T} \in \{0, -1\}^{(n - 2k) \times k}$ is the binary matrix that we use to control $\boldsymbol{P}$. By the forms of $\boldsymbol{P}$ and $\boldsymbol{A}$, given the problem instance $\boldsymbol{\pi}_{r, s}$, we have $\boldsymbol{AP} = \begin{bmatrix}
        \boldsymbol{I}_k \\
        -\boldsymbol{I}_k
    \end{bmatrix},$ and therefore the constraints $\boldsymbol{AP} \boldsymbol{y} \leq \boldsymbol{b}_s$ implies $\boldsymbol{y}_j = 0$ for $j = 1, \dots, k$ if $j \neq s$, and $\boldsymbol{y}_s \in [0, 1]$. Bedsides, the objective of the problem instance $\boldsymbol{\pi}_{r, s}$ is $\boldsymbol{y}^\top \boldsymbol{P}^\top \boldsymbol{Q}\boldsymbol{P}\boldsymbol{y}+  \boldsymbol{c}^\top \boldsymbol{P}\boldsymbol{y} =\boldsymbol{c}_r^\top \boldsymbol{P}\boldsymbol{y} =\boldsymbol{T}_{r, s}\boldsymbol{y}_s$, where $\boldsymbol{T}_{r, s}$ is the entry of matrix $\boldsymbol{T}$ in the $r^{th}$ row and $s^{th}$ column. Since $\boldsymbol{T}_{r, s} \in \{-1, 0\}$, and $\boldsymbol{y} \in [0, 1]$, we have the optimal objective of the PQP corresponding to the OQP $\boldsymbol{\pi}_{r, s}$ and the projection matrix $\boldsymbol{P}$ is $\boldsymbol{T}_{r, s}$. Therefore, for the set of QP problem instances $\{\boldsymbol{\pi}_{r, s}\}_{r \in \{1, \dots, n - 2k\}, s \in \{1, \dots, k\}}$, we choose the set of real-valued thresholds $\{\tau_{r, s}\}_{r \in \{1, \dots, n - 2k\}, s \in \{1, \dots, k\}}$, where $\tau_{r, s} = -\frac{1}{2}$. Then, for each subset $\boldsymbol{I} \subset \{1, \dots, n - 2k\} \times \{1, \dots,k\}$, we construct $\boldsymbol{P}$ by choosing $\boldsymbol{T}$ such that $\boldsymbol{T}_{r, s} = -1$ if $(r, s) \in \boldsymbol{I}$ and $\boldsymbol{T}_{r, s} = 0$ otherwise. Therefore
    \[
        \ell_{\boldsymbol{P}}(\boldsymbol{\pi}_{r, s}) \geq \tau_{r, s} \,\, \textup{if} \,\, (r, s) \in \boldsymbol{I},
        \] 
    and 
    \[
    \ell_{\boldsymbol{P}}(\boldsymbol{\pi}_{r, s}) < \tau_{r, s} \,\, \textup{otherwise}.
    \]
    This means that the function class can shatter the set of QP problem instances $\{\boldsymbol{\pi}_{r, s}\}_{r, s}$ above, and therefore $\Pdim(\cL) = \Omega(nk)$.  
\qed

\begin{algorithm}[tb] 
        \caption{The unrolled active set method $\Gamma_{\boldsymbol{\pi}}$ corresponding to the  OQP $\boldsymbol{\pi} = (\boldsymbol{Q}, \boldsymbol{c}, \boldsymbol{A}, \boldsymbol{b})$}
            \label{alg:unrolled-active-set-method-modified}
            \textbf{Input}: Projection matrix $\boldsymbol{P} \in \bbR^{n \times k}$\\
            \textbf{Output}: A recovered sub-optimal solution $\boldsymbol{P}\boldsymbol{y}_\textup{cand}$
             \begin{algorithmic}[1]
            \STATE Set $\tilde{\boldsymbol{Q}} = \boldsymbol{P}^\top \boldsymbol{Q} \boldsymbol{P}$, $\tilde{\boldsymbol{A}} = \boldsymbol{AP}$, and $\tilde{\boldsymbol{c}} = \boldsymbol{P}^\top \boldsymbol{c}$.
            
            \FOR{potential active set $\cA \subseteq \{1, \dots, m\}$, $\abs{\cA} \leq \min\{m, k\}$}
                \STATE Construct KKT matrix $\boldsymbol{K} = \begin{bmatrix}
                    \tilde{\boldsymbol{Q}} & \tilde{\boldsymbol{A}}_\cA^\top \\
                    \tilde{\boldsymbol{A}}_\cA & \boldsymbol{0}.
                \end{bmatrix}$
                
                \IF{$\det(\boldsymbol{K}) \neq 0$}
                    \STATE Compute $\begin{bmatrix}
                        \boldsymbol{y}_\textup{cand} \\
                        \boldsymbol{\lambda}_\textup{cand}
                    \end{bmatrix} = \boldsymbol{K}^{-1}\begin{bmatrix}
                        - \tilde{\boldsymbol{c}} \\
                        \boldsymbol{b}_\cA
                    \end{bmatrix}$.
    
                    \LineComment{Checking feasibility of potential solution $\boldsymbol{y}_\textup{cand}$}
                    \STATE $yFeasible = True$ 
                    \FOR{$j \not \in \cA$}
                        \IF{$\tilde{\boldsymbol{A}}_j^\top \boldsymbol{y}_\textup{cand} > \boldsymbol{b}_j$}
                            \STATE $yFeasible = False$
                            \STATE \textbf{break}
                        \ENDIF
                    \ENDFOR
                    \IF{$yFeasible$}
                        \LineComment{Checking validation of Lagrangian $\boldsymbol{\lambda}_{\textup{cand}}$}
                        $lambdaValid = True$
                        \FOR{$j \in \cA$}
                            \IF{$\boldsymbol{\lambda}_{\textup{cand}, j} < 0$}
                                \STATE $lambdaValid = False$
                                \STATE \textbf{break}
                            \ENDIF
                        \ENDFOR
                        
                        \IF{$lambdaValid$}
                            \STATE \textbf{return} $\boldsymbol{P}\boldsymbol{y}_\textup{cand}$
                        \ENDIF
                    \ENDIF
                \ENDIF
            \ENDFOR
        
        \end{algorithmic}
    \end{algorithm}

\section{Proofs for Section \ref{sec:extension}} \label{apx:extension}

\subsection{Proofs for Section \ref{sec:match-solution}} \label{apx:match-solution}
In this section, we will present the formal proof for Theorem \ref{thm:learning-to-match-solution}. First, note that under Assumption \ref{asmp:Q-PD}, given a QP problem instance $\boldsymbol{\pi} = (\boldsymbol{Q}, \boldsymbol{c}, \boldsymbol{A}, \boldsymbol{b})$, the objective matrix $\boldsymbol{Q}$ is already positive definite, ensuring that the optimal solution $\boldsymbol{\pi}^*$ is unique when combining with Assumption \ref{asmp:regularity}. 

Using the fact above, we will slightly modify the unrolled active set method (Algorithm \ref{alg:unrolled-active-set-method}) so that it corresponds to the OQP, instead of the perturbed OQP, takes input as a projection matrix $\boldsymbol{P}$ and output the recovered solution $\boldsymbol{P}\boldsymbol{y}^*$ from optimal solution $\boldsymbol{y}^*(\boldsymbol{P}, \boldsymbol{\pi})$ of the PQP corresponding to $\boldsymbol{P}$ and $\boldsymbol{\pi}$. The detailed modification is demonstrated in Algorithm \ref{alg:unrolled-active-set-method-modified}.

We first show that Algorithm \ref{alg:unrolled-active-set-method-modified} correctly output the optimal solution for the PQP corresponds to the OQP $\boldsymbol{\pi}$ and the projection matrix $\boldsymbol{P}$.
\begin{proposition} \label{prop:match-solution-correctness}
    Given a OQP $\boldsymbol{\pi}$, the algorithm $\Gamma_{\boldsymbol{\pi}}$, described by Algorithm \ref{alg:unrolled-active-set-method-modified} and corresponding to $\boldsymbol{\pi}$, correctly computes the optimal solution $\boldsymbol{y}^*(\boldsymbol{P}, \boldsymbol{\pi})$, i.e., $\boldsymbol{y}_\textup{cand} = \boldsymbol{y}^*(\boldsymbol{P}, \boldsymbol{\pi})$, and therefore the output $\boldsymbol{P}\boldsymbol{y}_\textup{cand}$ is the recovered solution.
\end{proposition}
\proof
    Again, the proof idea is similar to that of Lemma \ref{lm:unrolled-active-method-correctness}. For the \textbf{existence} part, from Lemma \ref{lm:localizing-solution-pertubred-PQP}, given the (unique) optimal solution of the PQP $\boldsymbol{y}^*(\boldsymbol{P}, \boldsymbol{\pi})$, there exists a subset $\cB \subseteq \cA(\boldsymbol{y}^*(\boldsymbol{P}, \boldsymbol{\pi}))$ of the active set corresponding to $\boldsymbol{y}^*(\boldsymbol{P}, \boldsymbol{\pi})$ such that $\tilde{\boldsymbol{A}}$ has linearly independent rows, and that $\boldsymbol{y}^*(\boldsymbol{P}, \boldsymbol{\pi})$ is the unique optimal solution of the new equality-constrained QP:
    \[
        \boldsymbol{y}^*(\boldsymbol{P}, \boldsymbol{Q}) = \min_{\boldsymbol{y} \in \bbR^k} \left\{\frac{1}{2}\boldsymbol{y}^*\tilde{\boldsymbol{Q}} \boldsymbol{y} + \tilde{\boldsymbol{c}}^\top \boldsymbol{y} \mid \tilde{\boldsymbol{A}}_\cB \boldsymbol{y} = \boldsymbol{b}_\cB\right\}.
    \]
    Note that Algorithm \ref{alg:unrolled-active-set-method-modified} considers all subsets $\cA \subset \{1, \dots, m\}$ that have at most $\abs{\cA} \leq k$ elements, and it will eventually check $\cB$. And when $\cA = \cB$, we can easily show that $(\boldsymbol{y}_\textup{cand}, \boldsymbol{\lambda}_\textup{cand})$ will pass all primal and dual feasibility checks, which means that $\boldsymbol{y}_\textup{cand} = \boldsymbol{y}^*(\boldsymbol{P}, \boldsymbol{\pi})$. Moreover, we can easily verify that $(\boldsymbol{y}_\textup{cand}, \boldsymbol{\lambda}_\textup{cand})$ is a KKT point of the equality-constrained QP above, meaning that $\boldsymbol{y}_\textup{cand}$ is also its unique optimal solution.

    For the \textbf{correctness} part, also similar to Lemma \ref{lm:unrolled-active-method-correctness}, we also show that given the output $(\boldsymbol{y}_\textup{cand}, \boldsymbol{\lambda}_\textup{cand})$, we can construct the point $(\boldsymbol{y}_\textup{cand}, \boldsymbol{\lambda})$, where $\boldsymbol{\lambda}_i = 0$ if $i \not \in \cA$ and $\boldsymbol{\lambda}_i = \boldsymbol{\lambda}_{\textup{cand}, i}$ if $i \in \cA$, that satisfies the KKT conditions of the PQP. This means that $\boldsymbol{y}_{\textup{cand}}$ is the optimal solution of the PQP. Besides, it is easy to check that $(\boldsymbol{y}, \boldsymbol{\lambda})$ also satisfies the KKT conditions of the equality-constrained QP corresponding with the same objective of PQP and the constraints $\tilde{\boldsymbol{A}}_\cA \boldsymbol{y} = \boldsymbol{b}_\cA$.
\qed

    Secondly, we will show that Algorithm \ref{alg:unrolled-active-set-method-modified} is also a GJ algorithm with bounded complexities. Again, the proof is similar to the proof of Lemma \ref{lm:gj-complexity-unrolled-active-set}.

    \begin{proposition} \label{prop:match-solution-complexity}
        Given a OQP $\boldsymbol{\pi}$, the algorithm $\Gamma_{\boldsymbol{\pi}}$ described by Algorithm \ref{alg:unrolled-active-set-method-modified} is a GJ algorithm with degree $\cO(m + k)$ and predicate complexity $\cO\left(m\min\left(2^m, (\frac{em}{k})^k\right) \right)$.
    \end{proposition}
    \proof
        Similar to the proof of Lemma \ref{lm:gj-complexity-unrolled-active-set}, we can claim that all the intermediate values computed by Algorithm \ref{alg:unrolled-active-set-method-modified} are all rational functions of (the entries of) $\boldsymbol{P}$ of degree $\cO(m + k)$. Moreover, we can also bound the number of distinct predicates (rational functions involved in the condition statements) by $\cO(m \cdot 2^m)$. 
    \qed

    Finally, we can formalize the proof of Theorem \ref{thm:learning-to-match-solution}.
    \proof[Proof of Theorem~\ref{thm:learning-to-match-solution}]
        This is a direct consequence from Proposition \ref{prop:match-solution-complexity} and Theorem \ref{thm:gj}.
    \qed 

\subsection{Proofs for Section \ref{sec:input-aware}}
We first formalize the following structural result, which says that given a set of $N$ input problem instances $\boldsymbol{\pi}_1, \dots, \boldsymbol{\pi}_N$, the outputs $f_{\boldsymbol{\theta}}(\boldsymbol{\pi}_i)$ admits piecewise polynomial structure, with bounded number of pieces. 
\begin{proposition} \label{lm:GJ-input-aware-f}
     Given any set of $N$ OQPs $\boldsymbol{\pi}_1, \dots, \boldsymbol{\pi}_N$, we can partition the space $\Theta \subset \bbR^{W}$ of neural network parameters into connected components $\{\cC_1, \dots, \cC_C\}$, where
     \[
        C \leq 2^{L + 1} \left(\frac{2eN(U + 2nk)}{W}\right)^{(L + 1)W}.
     \]
     Given a connected component $\cC_i$, the projection matrix $\boldsymbol{P}_{\boldsymbol{\pi}_i} = f_{\boldsymbol{\theta}}(\boldsymbol{\pi}_{i, \textup{flat}})$, for any $i \in \{1, \dots, N\}$, is a matrix with polynomials entries (in the neural network parameters $\boldsymbol{\theta}$) of degree at most $L + 2$.
\end{proposition}
\proof
    The result is a direct consequence of the result by \citet{anthony2009neural}, later adapted to the context of data-driven algorithm selection by \citet{cheng2024sample} (see e.g., Theorem 2.6). We simply adapt the results in the context of input-aware data-driven learning, the projection matrix, where the output of the neural network is $nk$, which is the size of the projection matrix. 
\qed

We now can present the proof of Theorem \ref{thm:input-aware-setting}.
\proof[Proof of Theorem \ref{thm:input-aware-setting}]
    Given $N$ OQPs $\boldsymbol{\pi}_1, \dots, \boldsymbol{\pi}_N$ and $N$ real-valued thresholds $\tau_1, \dots, \tau_N$, we first need to bound the number of sign pattern 
    \[
        \{\sign(\ell(f_{\boldsymbol{\theta}}(\boldsymbol{\pi}_{1, \textup{flat}})) - \tau_1), \dots, \sign(\ell(f_{\boldsymbol{\theta}}(\boldsymbol{\pi}_{N, \textup{flat}})) - \tau_N) \mid \boldsymbol{\theta} \in \Theta\}
    \]
    when varying $\boldsymbol{\theta} \in \Theta$.

    From Proposition \ref{lm:GJ-input-aware-f}, the parameter space $\Theta$ can be partitioned into connected components $\{\cC_1, \dots, \cC_C\}$, where
    \[
        C \leq 2^{L + 1} \left(\frac{2eN(U + 2nk)}{W}\right)^{(L + 1)W},
     \]
     and in each connected component $\boldsymbol{\cC} \subset \Theta$, the projection matrix $\boldsymbol{P}_{\boldsymbol{\pi}_i} = f_{\boldsymbol{\theta}}(\boldsymbol{\pi}_{i, \textup{flat}})$, for any $i \in \{1, \dots, N\}$, is a matrix with polynomials entries (in the neural network parameters $\boldsymbol{\theta}$) of degree at most $L + 2$. Now, in each connected components $\cC$, from Lemma \ref{lm:gj-complexity-unrolled-active-set}, the sign $\sign(\ell(f_{\boldsymbol{\theta}}(\boldsymbol{\pi}_{1, \textup{flat}}) - \tau_1)$ is determined by at most $m t $ polynomials, each of degree at most $\cO((L + 2)(m + k))$, where $t = \min\left(2^m, \left(\frac{emk}{k}\right)^k\right)$. Therefore, the number of signs 
     \[
        \{\sign(\ell(f_{\boldsymbol{\theta}}(\boldsymbol{\pi}_{1, \textup{flat}})) - \tau_1), \dots, \sign(\ell(f_{\boldsymbol{\theta}}(\boldsymbol{\pi}_{N, \textup{flat}})) - \tau_N) \mid \boldsymbol{\theta} \in \Theta\}
    \]
    acquired by varying $\boldsymbol{\theta} \in C$ is at most 
    \[
        \cO\left(\frac{8eNmt(L + 2)(m + k)}{W}\right)^W.
    \]
    This means that the number of signs 
    \[
        \{\sign(\ell(f_{\boldsymbol{\theta}}(\boldsymbol{\pi}_{1, \textup{flat}}) - \tau_1), \dots, \sign(\ell(f_{\boldsymbol{\theta}}(\boldsymbol{\pi}_{N, \textup{flat}}) - \tau_N) \mid \boldsymbol{\theta} \in \Theta\}
    \]
    acquired by varying $\boldsymbol{\theta} \in \Theta$ is at most 
    \[
        Z(N) = \cO\left(\frac{8eNmt(L + 2)(m + k)}{W}\right)^W \cdot 2^{L + 1} \left(\frac{2eN(U + 2nk)}{W}\right)^{(L + 1)W}.
    \]
    Solving the inequality $2^N \leq Z(N)$, and use the inequality $\log z \leq \frac{z}{\lambda} + \log\frac{\lambda}{e}$ for $z > 0$ and $\lambda > 0$ yields
    \[
        N = \cO(WL\log(U + mk) + W\min(m, k\log m)). 
    \]
    This completes the proof.
\qed

{
    \section{Gradient Update for Data-driven Learning the Projection Matrix for QPs}
    \label{apx:gradient-update}
    In this section, we will formalize derive the gradient update for learning the projection matrix for QPs in the data-driven framework. Recall that given a problem instance $\boldsymbol{\pi} = (\boldsymbol{Q}, \boldsymbol{c}, \boldsymbol{A}, \boldsymbol{b})$ and a projection matrix $\boldsymbol{P}$, we have
    \[
        \ell(\boldsymbol{P}, \boldsymbol{\pi)} = \min_{\boldsymbol{y} \in \bbR^k} \frac{1}{2}\boldsymbol{y}^\top \boldsymbol{P}^\top \boldsymbol{Q} \boldsymbol{P} \boldsymbol{y} + \boldsymbol{c}^\top \boldsymbol{P}\boldsymbol{y} \quad \text{s.t.} \quad \boldsymbol{APy} \leq \boldsymbol{b}.
    \]
    To calculate $\nabla_{\boldsymbol{P}}\ell(\boldsymbol{P}, \boldsymbol{\pi})$, we first recall the Envelope theorem. 
    \begin{lemma}[Envelope theorem, \cite{milgrom2002envelope}] \label{lm:envelope}
        Let $f(\boldsymbol{x}, \boldsymbol{\alpha})$ and $g_j(\boldsymbol{x}, \boldsymbol{\alpha)}$, where $j = 1, \dots, m$ be real-valued continuously differentiable function, where $\boldsymbol{x} \in \bbR^n$ and variables, and $\boldsymbol{\alpha} \in \bbR^l$ are parameters, and consider the parametric optimization problem
        \[
            \ell(\boldsymbol{\alpha}) = \min_{\boldsymbol{x}} f(\boldsymbol{x}, \boldsymbol{\alpha}) \quad \textup{subject to} \quad g_i(\boldsymbol{x}, \boldsymbol{\alpha)} \leq 0, i = 1, \dots, m.
        \]
        Let $\cL(\boldsymbol{x}, \boldsymbol{\alpha}, \boldsymbol{\lambda)}$ be the corresponding Lagrangian
        \[
            \cL(\boldsymbol{x}, \boldsymbol{\alpha}, \boldsymbol{\lambda}) = f(\boldsymbol{x}, \boldsymbol{\alpha}) + \sum_{i = 1}^m \lambda_i g_i(\boldsymbol{x}, \boldsymbol{\alpha}),
        \]
        where $\boldsymbol{\lambda}$ is the Lagrangian multiplier. Let $\boldsymbol{x}^*(\alpha)$, $\boldsymbol{\lambda}^*(\alpha)$ be the solution that minimizes the objective subject to the constraints, and let $\cL^*(\boldsymbol{\alpha}) = \cL(\boldsymbol{x}^*(\boldsymbol{\alpha}), \boldsymbol{\alpha}, \boldsymbol{x}^*(\boldsymbol{\alpha)})$.
        Assume that $\ell(\boldsymbol{\alpha})$ and $\cL^*(\boldsymbol{\alpha})$ are continuously differentiable, then
        \[
            \nabla_{\boldsymbol{\alpha}}\ell(\boldsymbol{\alpha)} =  \nabla_{\boldsymbol{\alpha}}\cL(\boldsymbol{x}, \boldsymbol{\alpha}, \boldsymbol{\lambda})|_{\boldsymbol{x} = \boldsymbol{x}^*(\boldsymbol{\alpha}), \boldsymbol{\lambda}^*(\boldsymbol{\alpha})}.
        \]
    \end{lemma}
    Assuming that the regularity condition holds, using Lemma \ref{lm:envelope}, we have
    \[
        \nabla_{\boldsymbol{P}}\ell(\boldsymbol{P}, \boldsymbol{\pi}) = (\boldsymbol{QP}\boldsymbol{y}^*(\boldsymbol{P}) + \boldsymbol{c} + \boldsymbol{A}^\top \boldsymbol{\lambda}^*(\boldsymbol{P}))\boldsymbol{y}^*(\boldsymbol{P})^\top. 
    \]

}

\end{document}